\documentclass[a4paper,twoside,10pt]{article}
\usepackage[a4paper,width=150mm,top=25mm,bottom=25mm,bindingoffset=6mm]{geometry}

\usepackage[english]{babel}

\usepackage{graphicx}
\usepackage[font=small,labelfont=bf]{caption}
\usepackage{subcaption}
\usepackage{color}
\usepackage{mdframed}
\usepackage{framed} 
\usepackage{cancel}
\usepackage[normalem]{ulem}

\usepackage{physics}
\usepackage{tabularx}
\usepackage{multirow}
\usepackage{makecell}
\usepackage{booktabs}
\usepackage{setspace}
\usepackage{lineno,hyperref}
\modulolinenumbers[5]

\usepackage{amsmath,amssymb,amsthm,mathtools,mathrsfs,epsfig,amsfonts,wasysym}
\usepackage{MnSymbol}

\usepackage[normalem]{ulem}
\usepackage{ragged2e}
\usepackage{csquotes}
\usepackage{tasks}
\usepackage{paralist}
\usepackage[inline]{enumitem}
\usepackage{tikz-cd}
\usepackage{stmaryrd}

\usepackage[titletoc,title]{appendix}

\usepackage{longtable}

\usepackage[titletoc,title]{appendix}
\usepackage{url,hyperref}
\hypersetup{colorlinks=true,linkcolor=blue,filecolor=magenta,urlcolor=cyan}
\usepackage{nicematrix}
\usepackage{tikz}
\usetikzlibrary{fit}

\usepackage[all]{xy}

\newtheorem*{example*}{Example}
\newtheorem*{reminder*}{Reminder}
\newtheorem*{remark*}{Remark}
\newtheorem*{note*}{Note}
\usepackage{authblk}
\usepackage{todonotes}
\usepackage{orcidlink}

\begin{document}
%%%%%%%%%%%%%%%%%%%%%%%%%%%%%%%%%%%%%%%%%%%%%%%%%%%%%%%%%%%%%%%%%%%%%%%%%%%%%%%%%%%%%%%%%%%%%%%%%%%%%%%%

\title{Mathematical modeling of tumor–immune system interactions: the effect of rituximab on breast cancer immune response}

\author[1,2]{Vasiliki Bitsouni\,\orcidlink{0000-0002-0684-0583}\thanks{\texttt{vbitsouni@math.uoa.gr}}}
\author[2]{Vasilis Tsilidis\,\orcidlink{0000-0001-5868-4984}\thanks{\texttt{vtsilidis@outlook.com}}}
\affil[1]{Department of Mathematics, National and Kapodistrian University of Athens, Panepistimioupolis, GR-15784 Athens, Greece}
\affil[2]{School of Science and Technology, Hellenic Open University, 18 Parodos Aristotelous Str., GR-26335 Patras, Greece}

\date{}

\maketitle

\begin{abstract}
tBregs are a newly discovered subcategory of B regulatory cells, which are generated by breast cancer, resulting in the increase of Tregs and therefore in the death of NK cells. In this study, we use a mathematical and computational approach to investigate the complex interactions between the aforementioned cells as well as CD8$^+$ T cells, CD4$^+$ T cells and B cells. Furthermore, we use data fitting to prove that the functional response regarding the lysis of breast cancer cells by NK cells has a ratio-dependent form. Additionally, we include in our model the concentration of rituximab - a monoclonal antibody that has been suggested as a potential breast cancer therapy - and test its effect, when the standard, as well as experimental dosages, are administered.
\end{abstract}

\numberwithin{equation}{section}
%%%%%%%%%%%%%%%%%%%%%%%%%%%%%%%%%%%%%%%%%%%%%%%%%%%%%%%%%%%%%%%%%%%%%%%%%%%%%%%%%%%%%%%
\section{Introduction}\label{intro}
B cells play an important role in antibody production, regulation of T cells and activation of CD4$^+$ T cells \cite{murphy2016janeway}. They are characterized by inhomogeneity and depending on their function are categorized in different classes, such as plasma cells which produce antibodies \cite{nutt2015generation} and B regulatory cells (Bregs) which regulate the function of other immune cells \cite{rosser2015regulatory}.
  
The relationship between B cells and cancer, even though is as important as the relationship of T cells and NK cells with cancer, which is generally more common in scientific research, has only recently started being studied  \cite{guo}.
In recent years, the discovery of tumor-infiltrating B cells has sparked new research regarding their role in cancer \cite{linnebacher2012tumor}. 
More specifically in breast cancer, the function of B cells seems to be very complex and is still debatable with different studies indicating them either as positive \cite{mahmoud2012prognostic,iglesia2014prognostic,xu2018prognostic} or negative \cite{mohammed2012relationship,mohammed2013relationship,miligy2017prognostic} mediators of the disease or remaining neutral \cite{west2011tumor,eiro2012impact,thompson2016immune}. Due to their big significance and rich interactions with breast cancer, tumor-infiltrating B cells have been characterized as a ``new hallmark of breast cancer" \cite{shen2018new}. 
  
Recent scientific publications \cite{olkhanud2009breast,olkhanud2011tumor,  biragyn2014generation}, discovered a sub-population of Bregs, named tumor-evoked Bregs (tBregs), which are being generated by the existence of breast cancer. tBregs in turn, cause an increase in Treg population by helping the differentiation of CD4$^+$ T cells to Tregs, which kill NK cells causing breast cancer to metastasize to the lungs.

The authors of \cite{olkhanud2011tumor} concluded that tBregs need to be controlled in order for breast cancer to regress, thus, suggesting the anti-CD20 monoclonal antibody rituximab as a potential cure for some types of  breast cancer. Rituximab targets the CD20 protein, which is mainly found on the surface of B cells, binding with it and triggering B cell death \cite{bosch2014drugsRituximab}. It is used to combat blood cancers such as leukaemia and lymphoma, as well as autoimmune diseases such as rheumatoid arthritis \cite{en:Rituximab}.

Studies regarding the effect of B cell depletion in cancer have been mixed. In \cite{kim2008b}, mice bearing lung cancer were depleted of B cells, through the use of an anti-CD20 antibody, which slowed tumor growth. Additionally, when active immunotherapy was used in conjunction with the anti-CD20 antibody, the authors observed increased anti-tumor effects and CD8$^+$ T cell levels. On the other hand, in \cite{aklilu2004depletion} the authors treated fifteen renal cell carcinoma and six melanoma patients with rituximab and IL-2 and found that B cell depletion produces no different results on IL-2 therapy. Moreover, in \cite{candolfi2011b}, B-cell-depleted mice bearing glioblastoma and wild-type mice treated with an anti-CD20 antibody bearing glioblastoma were given a treatment that induces tumor regression in 60\% of wild-type mice. The treatment completely failed in both classes of mice, as mice were unable to exhibit clonal expansion of anti-tumor T cells. Thus, the authors noted that B cells play the role of antigen presenting cells.

Mathematical models studying the role of B cells in cancer are scarce. The few published mathematical models that study B cells, mainly focus on their ability to produce antibodies \cite{ghosh2018mathematical, dhar2020numerical} or the relationship between mature B cells and progenitor B cells in B-cell acute lymphoblastic leukaemia \cite{leon2021car, nanda2013b}, while as far as we know, a mathematical model studying their regulatory activity does not exist. 

The goal of this study is to develop a mathematical framework within which we can investigate the complex interactions between breast cancer and the immune system, including B cells and tBregs, in order to get a better understanding of their functions, as well as investigate the efficacy of a potential B-cell-depletion breast cancer therapy through the administration of rituximab.
To this end, we derive a new mathematical model consisting of a system of coupled nonlinear ordinary differential equations.
In our model, we describe the interactions between breast cancer cells, NK cells, CD8$^+$ T cells, CD4$^+$ T cells, Tregs, B cells and tBregs, as well as the total concentration of rituximab administered to the organism. 

As far as CD8$^+$ T cells are concerned, they have been included in various mathematical models \cite{dePillis2005validated,dePillis2009}. The inclusion of CD8$^+$ T cells in our model will allow us to study their vital role of tumor-lysing, as well as their interactions with Tregs and non-Treg CD4$^+$ T cells.

Non-Treg CD4$^+$ T cells play a big role in anti-tumor immunity since they induce the proliferation of CD8$^+$ T cells, as well as NK cells through the production of IL-2. They are also activated by B cells. Furthermore, the process in which tBregs induce the proliferation of Tregs, relies on tBregs converting non-Treg CD4$^+$ T cells to Tregs. Even though non-Treg CD4$^+$ T cells have been studied in various mathematical models, the models mainly studied their ability to produce the cytokine IL-2 and not their interactions with other immune cells \cite{castiglione2007cancer,anderson2015qualitative,wei2017periodically,makhlouf2020mathematical},  whereas other models that exist in the literature either study their IL-2 production along with Treg generation \cite{robertson2012mathematical}, or their ability to induce the proliferation of effector cells \cite{dong2014mathematical}. Hence, the inclusion of non-Treg CD4$^+$ T cells in our model will allow us to study their rich interactions with other immune cells.

The layout of this study is as follows. In Section \ref{mathmodel}, we describe in detail the new mathematical model for tumor-immune interactions. In Section \ref{SectionLargeModelParameterEstimation}, we derive the model parameters. In Section \ref{numerics}, we study the dynamics of our model using numerical simulations. Finally, in Section \ref{conclusion}, we conclude with a summary and discussion of the results.

\section{Mathematical model}\label{mathmodel}
In this section, we develop a mathematical model in an attempt to study the interactions between breast cancer cells and the various immune cells, including B cells and tBregs, as well as the effect of rituximab on breast cancer progression. 

As the biochemical cascade of events linked with cancer growth and immune response are vastly complex, we note that there is no catholic agreement on those events. Therefore, we base our model on the following published scientific propositions:

\begin{enumerate}
    \item Breast cancer grows logistically in the absence of an immune response, as also discussed in Section \ref{SectionLargeModelParameterEstimationTumourcells}.
    \item Breast cancer promotes the proliferation of tBregs \cite{olkhanud2011tumor}.
    \item tBregs promote Treg generation, by converting them from non-Treg CD4$^+$ T cells \cite{olkhanud2011tumor}.
    \item Tregs kill NK cells, which causes lung metastasis \cite{olkhanud2009breast, pedroza2013}.
    \item Tregs aggressively suppress the proliferation of CD8$^+$ T cells and non-Treg CD4$^+$ T cells, when cocultured \cite{liyanage}.
    \item Tregs inhibit the cytotoxic activity of NK cells \cite{trzonkowski}.
    \item Both NK cells and CD8$^+$ T cells directly kill breast cancer cells \cite{murphy2016janeway,abbas2014cellular}.
    \item CD4$^+$ T cells improve the efficiency of CD8$^+$ T cells in killing cancer cells \cite{haabeth2014cd4+}.
    \item CD4$^+$ T cells are required for the generation of CD8$^+$ T cells \cite{lai2011roles, keene1982helper}.
    \item NK cells stimulate the proliferation of CD8$^+$ T cells \cite{assarsson}.
    \item CD4$^+$ T cells stimulate the proliferation of NK cells by producing the cytokine IL-2 \cite{meropol1998evaluation,antony2005cd4+ProduceIL2}.
    \item Rituximab only affects non-tBregs B cells, as tBregs express CD20 in low levels \cite{bodogai2013tBregsExpressCD20AtLowLevels}.
\end{enumerate}

\begin{figure}[ht]
 \makebox[\textwidth][c]{\includegraphics[width=.97\textwidth]{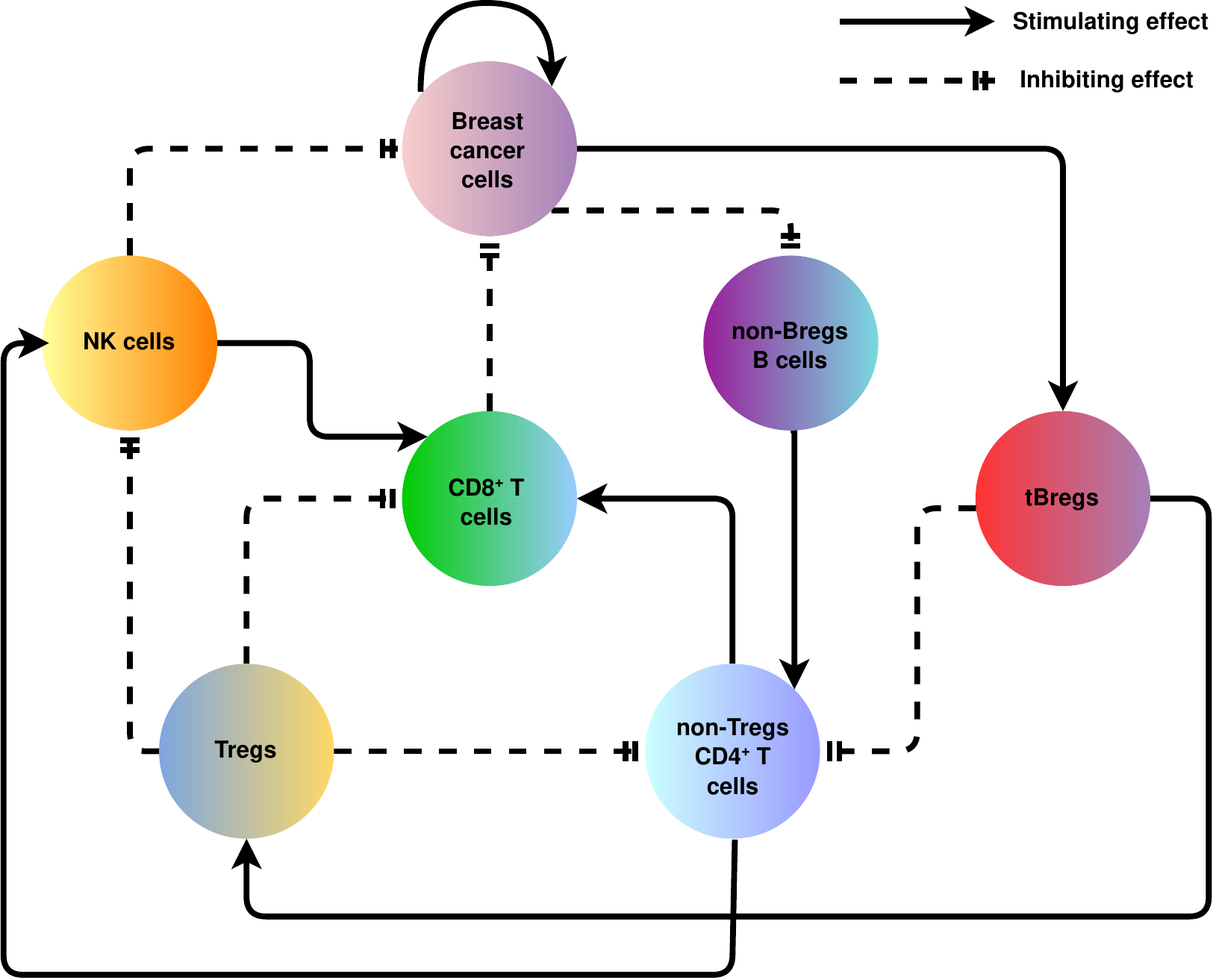}}%
 \caption{Interactions between the cells in system \eqref{LargeModelEquations}. Solid line: Stimulating effect. Dashed line: Inhibiting effect.}
 \label{fig:BigModel}
\end{figure}

Figure \ref{fig:BigModel} gives a schematic representation of the interactions between the cells in our model.

Therefore, our model considers seven cell populations and the concentration of rituximab. Let us denote by:
\begin{itemize}
    \item t, the time, measured in days.
    \item $T(t)$, the total tumor cell population, at time t.
    \item $N(t)$, the total NK cell population, at time t.
    \item $C(t)$, the total CD8$^+$ T cell population, at time t.
    \item $H(t)$, the total non-Treg CD4$^+$ T cell population, at time t.
    \item $R(t)$, the total Treg cell population, at time t.
    \item $B(t)$, the total non-tBreg B cell population, at time t.
    \item $B_T(t)$, the total tBreg cell population, at time t.
    \item $X(t)$, the mass of rituximab per blood volume, measured in micrograms per milliliter, at time t.
\end{itemize}

Our model consists of the following system of coupled nonlinear ordinary differential equations:

\begin{subequations} \label{LargeModelEquations}
\begin{align}  
& \dv{T}{t} = aT(1-bT) -c e^{-\lambda_R R} \frac{N^\delta}{s_N T^\delta + N^\delta} T - d \frac{C^l}{s_C T^l + C^l} T  \,, \label{LargeModeldT}\\
& \dv{N}{t} = \sigma_N - \theta_N N - pTN - \gamma_N R^{\delta_N} N   + \kappa H N  \,,            \label{LargeModeldN}        \\
& \dv{C}{t} = \sigma_C  - \theta_C C - qTC - \gamma_C RC + rNT + \frac{j_CT}{k_C+T}C + \frac{\eta_1 H}{\eta_2+H}C \,, \label{LargeModeldC} \\
& \dv{H}{t} = \sigma_H -\theta_H H + \frac{j_H T}{k_H + T} B H - c_1 H B_T \,, \label{LargeModeldH} \\
& \dv{R}{t} = \sigma_R -\theta_R R + c_1 H B_T                               \,, \label{LargeModeldR}              \\
& \dv{B}{t} = \sigma_B  -\theta_B B - c_2 T B - \gamma_B X^2 B \,, \label{LargeModeldB} \\
& \dv{B_T}{t} = - \theta_{B_T} B_T + c_2 T B \,, \label{LargeModeldBT} \\ 
& \dv{X}{t} = - \theta_X X + v(t) \,, \label{LargeModeldX}
\end{align}
\end{subequations} 
along with the initial condition:
\begin{multline}
\left(T\left(0\right),N\left(0\right),C\left(0\right),H\left(0\right),R\left(0\right),B\left(0\right),B_T\left(0\right),X\left(0\right)\right)= \\ \left(T^0,N^0,C^0,H^0,R^0,B^0,{B_T}^0,X^0\right)\in{\left[0,\infty\right)}^8.
    \label{Model_ICs}
\end{multline}

Table \ref{LargeModelTermDescription} lists all of the terms of system \eqref{LargeModelEquations}, along with a brief description. We subsequently give a detailed description of each model term. 

In \textbf{equation} \eqref{LargeModeldT}, the first term, $aT(1-bT)$, models the logistic growth of breast cancer. The last term, $- d \frac{C^l}{s_C T^l + C^l} T$, describes the CD8$^+$ T cells killing of breast cancer cells. This predation term is of a Hill form, with the ratio of CD8$^+$ T cells to breast cancer cells as the Hill function variable. It was firstly used in \cite{dePillis2005validated}, and subsequently in various other models \cite{dePillis2009,makhlouf2020mathematical,dePillis2013}. The lysis rate of tumor cells due to CD8$^+$ T cells seems to be a function of their ratio, thus this Hill form is able to capture that dynamic \cite{dePillis2005validated}. 

The second term, $-c e^{-\lambda_R R} \frac{N^\delta}{s_N T^\delta + N^\delta} T$, models the breast cancer lysis due to NK cells with Treg inhibition. We have used the function $e^{-\lambda_R R}$ to model the Treg inhibition of NK-induced breast cancer cell lysis, as it is positive and it does not affect the lytic activity of NK cells when there are no Tregs. It has also been used in \cite{dePillis2013} for the same reasons. The Hill function $c \frac{N^\delta}{s_N T^\delta + N^\delta}$ is used in the same way as $d \frac{C^l}{s_C T^l + C^l}$ is used to model the NK-induced breast cancer cell lysis. The simpler functional response function $cN$ is used in various other models to capture the same dynamic \cite{dePillis2005validated,dePillis2009}. Nevertheless, data from \cite{shenouda2017exVivo} show us that this Hill term seems to be able to capture the NK-induced breast cancer lysis more accurately. For a more in-depth discussion, see Section \ref{SectionLargeModelParameterEstimationNKcells}. Modeling the lytic activity of NK cells using a Hill function is a novel approach, since as far as we know there does not exist a model using this rational Hill form for this purpose.

In \textbf{equation} \eqref{LargeModeldN}, the first term, $\sigma_N$, represents the constant source of NK cells from the organism, whereas the second term, $-\theta_N N$, represents the natural NK cell death. The third term, $-pTN$, represents the inactivation of NK cells after interacting with tumor cells. A similar inactivation term has been used in \cite{kuznetsov1994} for the case of CD8$^+$ T cells and in other models for the NK cells case such as in \cite{depillis2003mathematical} and \cite{makhlouf2020mathematical}. The forth term, $- \gamma_N R^{\delta_N} N$, is used to model the Treg-induced NK apoptosis. The form of this term is derived from data fitting experiments based on data from \cite{shenouda2017exVivo}. For a more in-depth discussion, see Section \ref{SectionLargeModelParameterEstimationNKcells}. The last term, $\kappa HN$, is used to model the fact that NK cells proliferate in the presence of the cytokine IL-2 \cite{meropol1998evaluation}. Since non-Treg CD4$^+$ T cells are the main producers of IL-2 \cite{antony2005cd4+ProduceIL2}, we use this term as a proxy due to our model not including IL-2.

In \textbf{equation} \eqref{LargeModeldC}, the first term, $\sigma_C$, represents the constant source of CD8$^+$ T cells from the organism, whereas the second term, $-\theta_C C$, represents the natural CD8$^+$ T cell death. The third term, $-qTC$, represents the inactivation of CD8$^+$ T cells due to their interaction with breast cancer cells. It has been used in various models, such as in \cite{kuznetsov1994} and \cite{depillis2003mathematical}. The forth term, $-\gamma_CRC$, is used to model the suppression of CD8$^+$ T cell proliferation by Tregs. In \cite{liyanage}, the authors found that when coculturing CD8$^+$ T cells with Tregs taken from pancreatic and breast cancer patients, Tregs suppressed the proliferation of CD8$^+$ T cells. The fifth term, $r N T$, represents CD8$^+$ T cell recruitment due to the debris from tumor cells lysed by NK cells \cite{assarsson,huang1994roleNKdebris} and has been used in various models \cite{dePillis2005validated,makhlouf2020mathematical}. The sixth term, $\frac{j_CT}{k_C+T}C$, models the activation of CD8$^+$ T cells due to the presence of breast cancer cells and is included since CD8$^+$ T cells are part of the adaptive immune system. It has the same form as in \cite{kuznetsov1994}. The final term, $\frac{\eta_1 H}{\eta_2+H}C$, represents the CD4$^+$ T-cell-induced CD8$^+$ T cell proliferation and it has a Michaelis-Menten form. CD4$^+$ T cells can directly help the activation of CD8$^+$ T cells through cell-cell interactions via the CD40-CD154 signal pathway or indirectly through the production of IL-2 \cite{lai2011roles}. 

In \textbf{equation} \eqref{LargeModeldH}, the first term, $\sigma_H$, represents the constant source of non-Treg CD4$^+$ T cells from the organism, whereas the second term, $-\theta_H H$, represents the natural non-Treg CD4$^+$ T cell death. The third term, $\frac{j_H T}{k_H + T} B H$, represents the proliferation of non-Treg CD4$^+$ T cells due to the existence of breast cancer and is included since CD8$^+$ T cells are part of the adaptive immune system. Non-tBreg B cells appear in this term as they activate non-Treg CD4$^+$ T cells, acting as antigen presenting cells and thus without them, non-Treg CD4$^+$ T cells would not get activated. The last term, $-c_1 H B_T$, represents the differentiation of non-Treg CD4$^+$ T cells to Tregs due to tBregs \cite{olkhanud2011tumor}. As non-Treg CD4$^+$ T cells are converted to Tregs and Tregs do not seem to play a part in this conversion, we choose to only include the non-Treg CD4$^+$ T cells and tBregs in this term. In equation \eqref{LargeModeldH}, this term has a negative sign since non-Treg CD4$^+$ T cells are decreasing during this procedure.

In \textbf{equation} \eqref{LargeModeldR}, the first term, $\sigma_R$, represents the constant source of Tregs from the organism, whereas the second term, $-\theta_R R$, represents the natural Treg death. The last term, $c_1 H B_T$, represents the conversion of non-Treg CD4$^+$ T cells to Tregs. This term is the opposite of equation's \eqref{LargeModeldH} corresponding conversion term, since Tregs are increasing during this procedure and we want the two terms to have the same absolute value, because the same number of non-Treg CD4$^+$ T cells that are lost, become Tregs. 

In \textbf{equation} \eqref{LargeModeldB}, the first term, $\sigma_B$, represents the constant source of non-tBreg B cells from the organism, whereas the second term, $-\theta_B B$, represents the natural non-tBreg B cells death. The third term, $-c_2 T B$, represents the breast-cancer-induced differentiation of non-tBreg B cells to tBregs \cite{olkhanud2011tumor}. Just like the conversion term $c_1 H B_T$ in equations \eqref{LargeModeldH} and \eqref{LargeModeldR}, we only include non-tBreg B cells and breast cancer cells in this term as only these two seem to play a role in the conversion. In equation \eqref{LargeModeldB}, this term has a negative sign since non-tBreg B cells are decreasing during this procedure. The last term, $-\gamma_B X^2 B$, represents the rituximab-induced non-tBreg B cell apoptosis. The trophic function of this term is chosen to be of power form, since the term gets zeroed when there is no rituximab in the organism and because it makes a good fit to data found in \cite{tobinai1998feasibilityRituximabHalfLife} and \cite{cooper2004efficacyRituximab}.  

In \textbf{equation} \eqref{LargeModeldBT}, there is no intrinsic growth term, since we assume that tBregs do not exist in the organism in the absence of breast cancer. The first term, $-\theta_{B_T} B_T$, represents the natural tBreg cell death. The last term, $c_2 T B$, represents the differentiation of non-tBregs B cells to tBregs. This term is the opposite of equation's \eqref{LargeModeldB} corresponding conversion term since tBregs are increasing during this procedure and we want the two terms to have the same absolute value, since the same number of non-tBreg B cells that are lost, become tBregs. 

In \textbf{equation} \eqref{LargeModeldX}, the first term, $-\theta_{X} X$, represents the excretion of rituximab from patients. The last term, $v(t)$, is a function of time that models the mass of rituximab per liter of blood that gets infused into a patient per amount of time and is measured in $\frac{\mu\text{g}}{\text{mL} \cdot \text{day}}$.

\begin{table}
\centering
\caption{Description of the terms of system \eqref{LargeModelEquations}.}
\label{LargeModelTermDescription}
\noindent\makebox[\textwidth]{
\begin{tabular}{ccm{12cm}} 
\toprule
\textbf{Deriv.}               & \textbf{Term}                                                  & \textbf{Description}                                                     \\ 
\midrule
\multirow{3}{*}{$\dv{T}{t}$}   & $aT(1-bT)$                                                      & Logistic tumor growth                                                    \\
                               & $- c e^{-\lambda_R R} \frac{N^\delta}{s_NT^\delta + N^\delta}T$ & NK-induced tumor death with Treg inhibition                              \\
                               & $- d \frac{C^l}{s_C T^l + C^l} T$                               & CD8$^+$ T-induced tumor death                                            \\ 
\midrule
\multirow{5}{*}{$\dv{N}{T}$}   & $\sigma_N$                                                      & Constant source of NK cells                                               \\
                               & $- \theta_N N$                                                  & Programmed NK cell death                                                  \\
                               & $-pTN$                                                          & NK death by exhaustion of tumor-killing resources                \\
                               & $\gamma_N R^{\delta_N} N$                                       & Treg-induced NK apoptosis                                                 \\
                               & $\kappa H N$                                                    & CD4$^+$ T-induced NK cell proliferation                                    \\ 
\midrule
\multirow{7}{*}{$\dv{C}{t}$}   & $\sigma_C$                                                      & Constant source of CD8$^+$T cells                                         \\
                               & $- \theta_C C$                                                  & Programmed CD8$^+$T cell death                                            \\
                               & $- qTC$                                                         & CD8$^+$T cells death from exhaustion of tumor-killing resources          \\
                               & $- \gamma_C RC$                                                 & Suppression of  the proliferation of CD8$^+$ T cells by Tregs                     \\
                               & $rNT$                                                           & CD8$^+$ T cell recruitment due to NK-lysed tumor debris                  \\
                               & $\frac{j_CT}{k_C+T}C $                                          & Activation of CD8$^+$ T cells due to the presence of breast cancer cells  \\
                               & $\frac{\eta_1 H}{\eta_2+H}C $                                   & CD4$^+$ T-induced CD8$^+$ T proliferation                                  \\ 
\midrule
\multirow{4}{*}{$\dv{H}{t}$}   & $\sigma_H$                                                      & Constant source of non-Treg CD4$^+$ T cells                                \\
                               & $-\theta_H H $                                                  & Programmed non-Treg CD4$^+$ T cell death                                   \\
                               & $\frac{j_H T}{k_H + T} B H$                                     & CD4$^+$ T cell recruitment due to breast cancer, with B cell help          \\
                               & $- c_1 H B_T $                                                  & Differentiation of non-Treg CD4$^+$ T cells to Tregs due to tBregs       \\ 
\midrule
\multirow{3}{*}{$\dv{R}{t}$}   & $\sigma_R $                                                     & Constant source of Tregs                                                  \\
                               & $-\theta_R R$                                                   & Programmed Treg death                                                     \\
                               & $c_1 H B_T $                                                    & Differentiation of non-Treg CD4$^+$ T cells to Tregs due to tBregs       \\ 
\midrule
\multirow{4}{*}{$\dv{B}{t}$}   & $\sigma_B$                                                      & Constant source of non-tBreg B cells                                      \\
                               & $-\theta_B B$                                                   & Programmed non-tBreg B cell death                                         \\
                               & $- c_2 T B $                                                    & Breast-cancer-induced differentiation of non-tBreg B cells to tBregs      \\
                               & $- \gamma_B X^2 B$                                              & Rituximab-induced non-tBreg B cell apoptosis                              \\ 
\midrule
\multirow{2}{*}{$\dv{B_T}{t}$} & $ - \theta_{B_T} B_T$                                           & Programmed tBreg death                                                    \\
                               & $c_2 T B$                                                       & Breast-cancer-induced differentiation of non-tBreg B cells to tBregs      \\ 
\midrule
\multirow{2}{*}{$\dv{X}{t}$}   & $- \theta_X X$                                                  & Excretion of rituximab                                                    \\
                               & $v(t)$                                                          & Rituximab injection                                                       \\
\bottomrule
\end{tabular}}
\end{table}

\section{Parameter estimation} \label{SectionLargeModelParameterEstimation}
In this section, we carefully determine the model parameters. Since the model features a large amount of parameters, we use different methods in order to determine them, such as finding their value in biological literature, data fitting them based on biological research, calculating them based on the biological homeostasis states we found in  Appendix \ref{SectionEquilibriumStates}, borrowing them from other mathematical models or estimating them in order for our model to exhibit biological reasonable results. Below we give a detailed explanation about each parameter. A summary of the description and value of each parameter of the model can be found in Table \ref{tab:2}.

\subsection{The tumor} \label{SectionLargeModelParameterEstimationTumourcells} 
The breast cancer growth rate, $a=0.17$ day$^{-1}$, and inverse of carrying capacity, $b=10^{-10}$ cell$^{-1}$, are found using Mathematica's \verb|NonlinearModelFit| function to fit the logistic growth equation to breast tumor growth data from NSG mice found in \cite{puchalapalli2016nsg}. The authors of \cite{puchalapalli2016nsg} compared the growth and metastasis of three different breast cancer cells lines, namely CN34BrM, MDA-231 and SUM1315 cell lines, on athymic nude mice and NSG mice. The difference between the two kinds of mice is that the former lack T cells, since they are athymic, while their innate immunity is intact meaning that they still have NK cells. On the contrary, NSG mice not only are depleted of T cells, but also of B cells, while their NK activity is extremely low, therefore having impaired innate immunity. As can be seen in  Figure 1 in \cite{puchalapalli2016nsg}, the study showed that tumor growth in NSG mice was greater compared to athymic nude mice, therefore providing a better model of breast cancer growth in an immunodeficient organism. Even though similar data fitting experiments have already been conducted \cite{sarapata2014comparison}, we choose to run our own data fitting experiments due to the superiority of NSG mice versus athymic mice or BALB/c mice - which have the same immune cells as athymic nude mice \cite{charlesriverlabs} - like the ones in which the authors of \cite{sarapata2014comparison} have based their data fitting experiments on.
    
\begin{figure}[t]
 \makebox[\textwidth][c]{\includegraphics[width=1\textwidth]{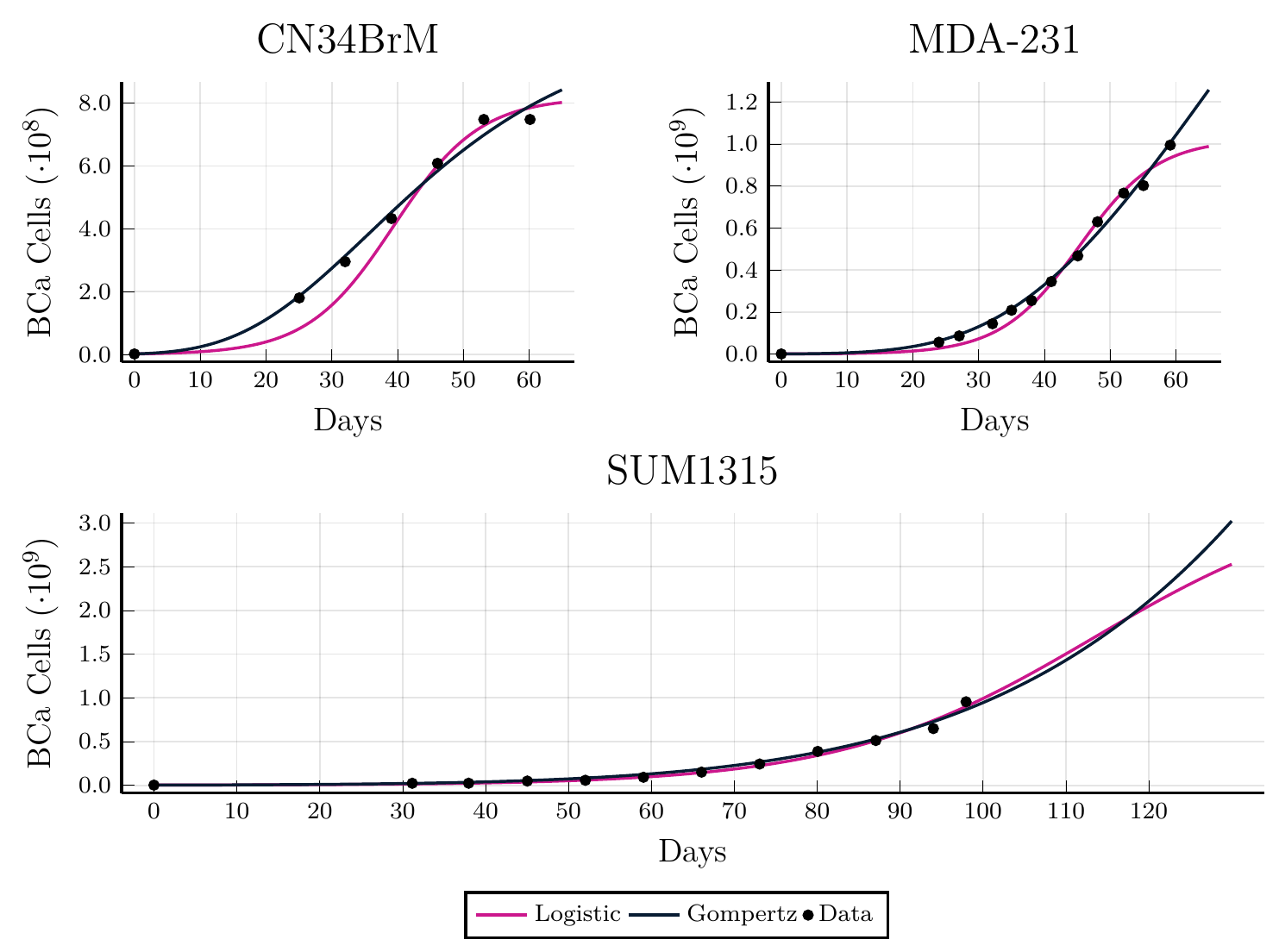}}%
  \caption{Fitting of the logistic and Gompertzian growth models to data from \cite{puchalapalli2016nsg}.}
  \label{fig:tumorGrowthFit}
\end{figure}
    
    Two of the most common mathematical models used to capture the growth of cancer cells are the logistic model 
 \begin{equation} \label{logisticEquation}
     \dv{p}{t} = r p \left(1-\frac{p}{K}\right)\,, \quad p(0) = p_0\,,
 \end{equation}
    and the Gompertzian model \cite{dePillis2014modelingBOOK}
    \begin{equation}
        \dv{p}{t} = rp \ln{\frac{K}{p}}\,, \quad \; p(0)=p_0\,,
    \end{equation}
    where $r$ is the intrinsic growth of the population and $K$ is its carrying capacity.

    Before we fit the two models to the data, we need to convert the data from tumor volume to total number of cancer cells that exist in each tumor volume. We use the same method as in \cite{dePillis2013}. Assuming a spherical tumor cell diameter of approximately 15.15$\mu$m yields a tumor cell volume of approximately $ 1.82 \cdot 10^3 \mu m^3$. Converting the data points from Figure 1 in \cite{puchalapalli2016nsg}, from mm$^3$ to $\mu$m$^3$ and dividing them by the tumor cell volume calculated, yields the total cancer cell number in each tumor volume.
        
    The results from fitting the above two growth models to the data from \cite{puchalapalli2016nsg} for the three breast cancer cell lines, are summarized in Figure \ref{fig:tumorGrowthFit} and Table \ref{tableLogisticGompertzianParameters}. As we can see in Figure \ref{fig:tumorGrowthFit}, both models make a very good fit to the data, which is consistent with the results from \cite{sarapata2014comparison}. With that in mind and considering that the logistic model is simpler and would make analysis easier, we pick the logistic over the Gompertzian function for our model.

\begin{table}
\centering
\caption{Parameter values for the logistic and Gompertzian model.}
\label{tableLogisticGompertzianParameters}
\begin{tabular}{ccccc} 
\toprule
\multirow{2}{*}{\textbf{Cell Line}} & \multicolumn{2}{c}{\textbf{Logistic Model}} & \multicolumn{2}{c}{\textbf{Gompertzian Model}}  \\ 
\cmidrule{2-5}
                                    & Growth Rate & Carrying Cap.                 & Growth Rate & Carrying Cap.                     \\ 
\midrule
CN34BrM                             & 0.16511     & $7.58 \cdot 10^8$             & 0.0513      & $1.05 \cdot 10^9$                 \\
MDA-231                             & 0.16835     & $1.03 \cdot 10^9$             & 0.0328      & $3.6 \cdot 10^9$                  \\
SUM1315                             & 0.06554     & $3.39 \cdot 10^9$             & 0.007       & $4.92 \cdot 10^{11}$                \\
\bottomrule
\end{tabular}
\end{table}

Finally, we choose the value of breast cancer growth rate to be the round up maximum value found by our data fitting experiments, that is 0.16835 day$^{-1}$. On the other hand, we chose the inverse of the carrying capacity to be a bit smaller than the lowest value found on our data fitting experiments which is approximately $2.8 \cdot 10^{-10}$ cell$^{-1}$. Our reasoning for doing so is because we want to study the case of an aggressive breast cancer. 

The value range for the maximum rate at which NK cells lyse cancer cells, $c$, the value of $(\frac{N}{T})^\delta$ for half-maximal NK toxicity, $s_N$, and the Hill coefficient, $\delta$, are determined through data fitting experiments based on data from \cite{shenouda2017exVivo}. In that study, the authors collected blood samples from normal donors and female breast cancer patients. The NK cells collected from their blood were expanded and subsequently placed on wells, along with either $2 \cdot 10^5$ cells of the triple negative breast cancer cell line MDA-MB-231 or $4 \cdot 10^5$ cells of the HER2-positive breast cancer cell line MDA-MB-453, at various ratios. After 4-5 hours the percent-specific lysis of breast cancer cells by NK cells was calculated. Assuming, that both cell populations are not able to grow inside the wells due to the lack of nutrients and space, we can use the following initial value problem to model the described phenomenon:

\begin{subequations} \label{tumorLysisByNKSystemDataFit}
\begin{align}
    \dv{T}{t} &= - f\left(T,N\right) T\left(t\right)\,, \quad T(0) = T_E\,, \\
    \dv{N}{t} &= -\theta_{N_E} N(t)\,, \quad N(0) = ratio \cdot T_E\,,
\end{align}
\end{subequations}
where $T$ is the breast cancer cell population, $N$ is the NK cell population, $\theta_{N_E}$ is the rate of natural NK cell death \textit{in vitro}, $T_E$ is the initial number of breast cancer cells, $ratio$ is the ratio of NK cells to breast cancer cells and $f\left(T,N\right)$ is the trophic function describing the killing of breast cancer cells by NK cells.

Using data from Figure 5 in \cite{olkhanud2009breast} we get that $\theta_{N_E} = 0.7414$ day$^{-1}$ (see Appendix \ref{app:naturaldeath} for more on how to calculate the turnover rate \textit{in vitro}). 

As far as the trophic function $f\left(T,N\right)$ is concerned, we use three different functions in order to determine which one makes the best fit and is therefore able to capture the dynamics of NK cells killing breast cancer cells more accurately. We use a power form, a rational Hill form and a Michaelis-Menten form, thus during our data fitting experiments the trophic function $f\left(T,N\right)$ takes one of the following forms:

\begin{equation} \label{trophicFunctionsRelation}
    f\left(T,N\right) = c N^\delta \text{ \quad or \quad} = c\frac{N^\delta}{s_N T^\delta + N^\delta} \text{ \quad or \quad}   = c \frac{N}{\delta + N}\,.
\end{equation}

Using Mathematica's \verb|ParametricNDSolveValue| function, we are able to solve problem \eqref{tumorLysisByNKSystemDataFit} numerically and get the percent-specific lysis of breast cancer cells by NK cells as a function of the form

\begin{equation} \label{PercentSpecificLysisRelation}
    \frac{1-T\left(t_{\textit{final}}\right)}{T_E}\,,
\end{equation}
where $t_{\textit{final}} = 5/24$ days, since we assume that authors kept the cells in the wells for 5 hours. Function \eqref{PercentSpecificLysisRelation} depends on the parameters of problem \eqref{tumorLysisByNKSystemDataFit}, as well as the $ratio$. We then pass function \eqref{PercentSpecificLysisRelation} on to Mathematica's \verb|NonlinearModelFit| function, which allows us to fit function \eqref{PercentSpecificLysisRelation} to the percent-specific lysis of breast cancer cells by NK cells expanded from breast cancer patients data taken from Figure 2 in \cite{shenouda2017exVivo}, with $ratio$ as the independent variable.

\begin{table}[]
\centering
\caption{Parameter values admitted from data fitting  problem \eqref{tumorLysisByNKSystemDataFit} to data from \cite{shenouda2017exVivo}.}
\label{tumorLysisByNKParameterTable}
\begin{tabular}{lcc}
\toprule
                     & \textbf{MDA-MB-231/luc}  & \textbf{MDA-MB-453}     \\ \midrule
\multicolumn{1}{c}{} & \multicolumn{2}{c}{\textbf{Power Form}}            \\ \midrule
$c$                  & $1.462 \cdot 10^{-7}$    & $2.07 \cdot 10^{-5}$    \\
$\delta$             & 1.2089                   & 0.7883                  \\ \midrule
\multicolumn{1}{c}{} & \multicolumn{2}{c}{\textbf{Rational Form}}         \\ \midrule
$c$                  & 11.2263                  & 19.6448                 \\
$\delta$             & 1.33332                  & 0.8249                  \\
$s_N$                & 39.222                   & 3.85119                 \\ \midrule
\multicolumn{1}{c}{} & \multicolumn{2}{c}{\textbf{Michaelis-Menten Form}} \\ \midrule
$c$                  & 55.0679                  & 22.858                  \\
$\delta$             & $1.8547 \cdot 10^{7}$    & $1.23545 \cdot 10^{6}$  \\ \bottomrule
\end{tabular}
\end{table}

\begin{figure}[t]
 \makebox[\textwidth][c]{\includegraphics[width=1\textwidth]{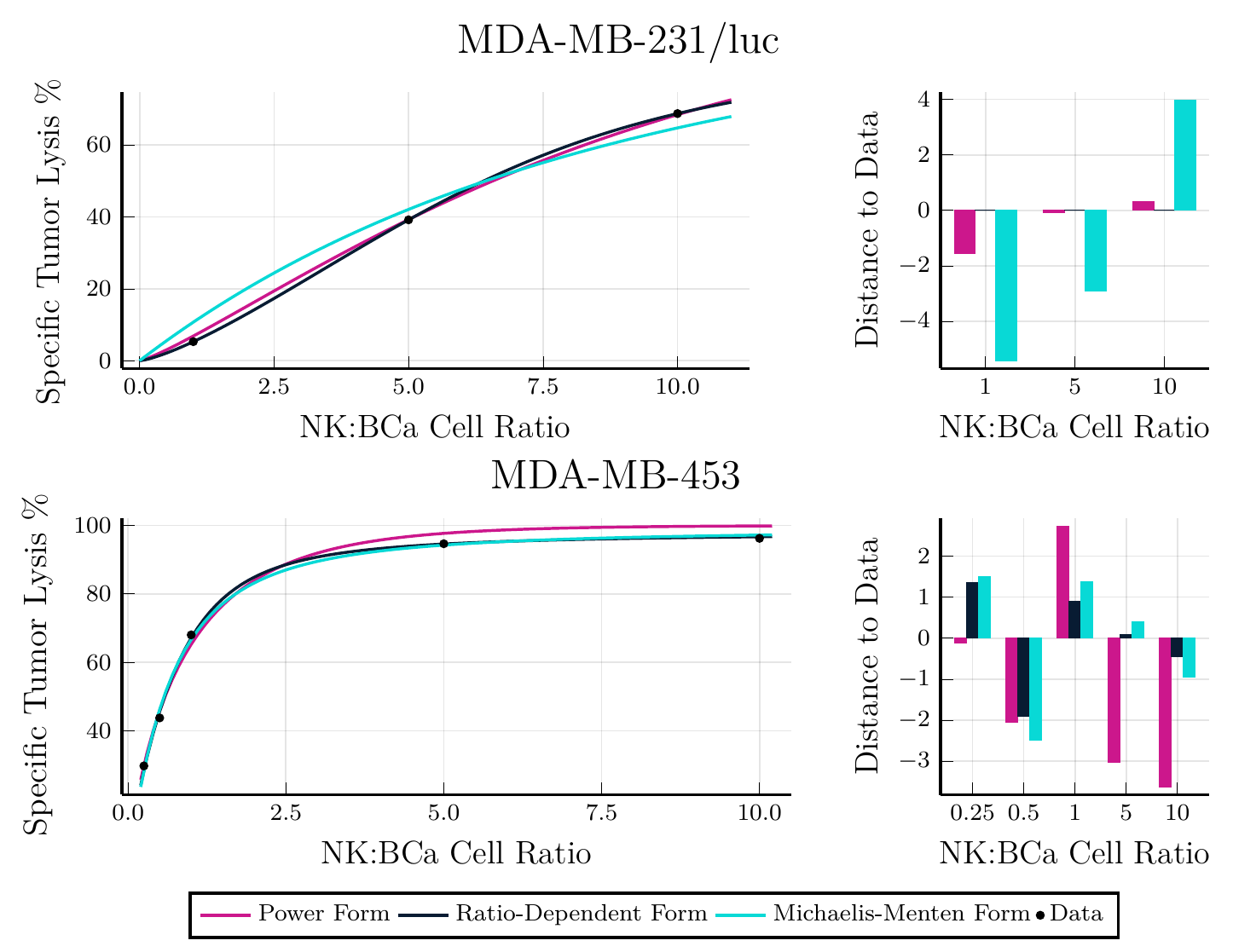}}%
  \caption{Left: Lysis curves of breast cancer cells by NK cells predicted by fitting the parameters of problem \eqref{tumorLysisByNKSystemDataFit} to data from \cite{shenouda2017exVivo}. Right: The residuals at each data point.}
  \label{fig:tumorLysisByNK}
\end{figure}

The percent-specific lysis curves predicted by  problem \eqref{tumorLysisByNKSystemDataFit} along with the distance to data at each data point are given in Figure \ref{fig:tumorLysisByNK}. Table \ref{tumorLysisByNKParameterTable} lists the parameters determined from our data fitting experiments. As we can see in Figure \ref{fig:tumorLysisByNK}, the rational Hill form makes the best fit regarding both breast cancer cell lines. Although the rational Hill form has one more variable and is therefore easier to be fitted, both breast cancer cell lysis curves seem to exhibit a saturation effect as the ratio gets larger, hence it is natural for the rational Hill form to make a better fit.

In \cite{dePillis2005validated}, a similar data fitting experiment was conducted, in which it was found that the power form, and more specifically a linear function, makes a very good fit in the case of NK cells lysing YAC-1 tumor cells. Different forms of the functional response function were not tested. This contrast between \cite{dePillis2005validated} and our simulations could be explained in two ways. Firstly, YAC-1 is a lymphoma cell line, unlike the two breast cancer cell lines we used in our simulations. Thus, it seems that NK cells lyse different cancer cell types in different ways. Secondly, as we already discussed, the cell lysis data we used in our simulations seem to exhibit a saturation effect as the ratio of NK to breast cancer cells gets larger, something that is not true with the respective data used in \cite{dePillis2005validated}, which also explains the  difference between the outcomes of \cite{dePillis2005validated} and our simulations. 

In our numerical simulations we vary those three parameters in order to study their effect on the breast cancer-immune dynamics.

The Treg-induced NK cell inhibition coefficient, $\lambda_R = 1 \cdot 10^{-8}$ cell$^{-1}$, is found to give the best fit to known data.

The maximum rate at which CD8$^+$ T cells lyse cancer cells, $d = 1.7$ day$^{-1}$, the value of $(\frac{C}{T})^l$ for half-maximal CD8$^+$ T cell toxicity, $s_C = 3.5 \cdot 10^{-2}$, and the Hill coefficient, $l = 1.7$, are borrowed from \cite{dePillis2013}, in which the authors derived the value of these parameters using data found in \cite{diefenbach} and \cite{dudley2002cancerdePillisData}.

\subsection{The NK cells} \label{SectionLargeModelParameterEstimationNKcells}
The constant source of NK cells, $\sigma_N = 1.13 \cdot 10^8$ cells $\cdot$ day$^{-1}$, is taken from \cite{zhang}. In that study, the authors found that healthy young adults have a total NK production rate of $(15 \pm 7.6) \cdot 10^6$ cells $ \cdot $ L$^{-1}$ · day$^{-1}$, while healthy older adults have one of $(7.3 \pm 3.7) \cdot 10^6$ cells $ \cdot $ L$^{-1}$ · day$^{-1}$. Considering that the average amount of blood in the human body is about 5 liters \cite{starr2012biology} and choosing the maximum NK production rate, we get the value of $\sigma_N$.

The rate of programmable NK cell death, $\theta_N = 0.06301 \text{ day}^{-1} $ is found by assuming the exponential decay of NK cells. Furthermore, the half-life of NK cells in humans is 1 to 2 weeks \cite{zhang}. Here, we choose an NK cell half-life of 11 days with a corresponding programmable NK death rate of $
    \theta_N = \frac{\ln2}{11 \text{ day}} \approx 0.06301 \text{ day}^{-1}\,. $

\begin{figure}
 \makebox[\textwidth][c]{\includegraphics[width=1\textwidth]{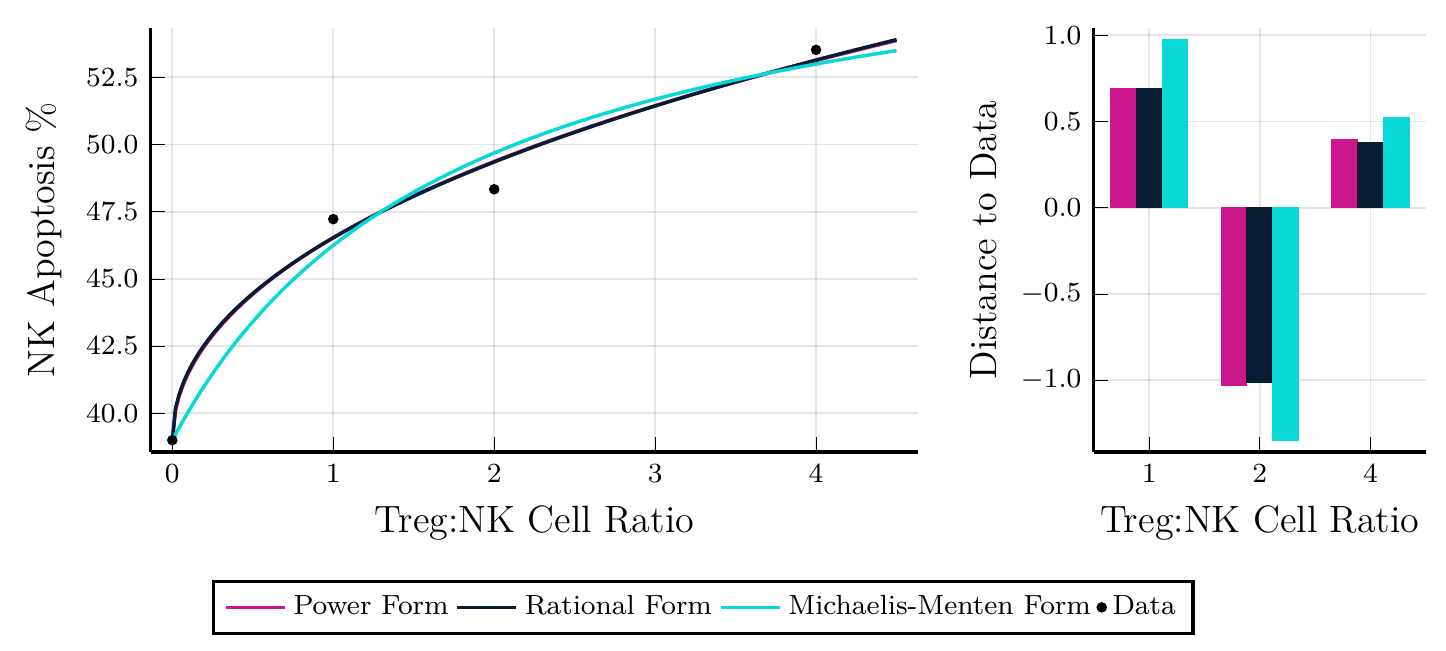}}%
  \caption{Left: Lysis curves of NK cells by Tregs predicted by fitting the parameters of problem \eqref{NKinducedApoptosisByTregsSystemDataFit} to data from \cite{olkhanud2009breast}. Right: The residuals at each data point.}
  \label{fig:NKlysisByTreg}
\end{figure}

The rate of Treg-induced NK cell apoptosis, $\gamma_N = 10^{-6}$ cell$^{-\delta_N}$ $ \cdot $ day$^{-1}$, and NK cell resistance to Treg-induced apoptosis coefficient, $\delta_N=0.5$, is determined by fitting data found in \cite{olkhanud2009breast}. In particular, the authors of \cite{olkhanud2009breast} cultured NK cells along with Tregs in wells, at various ratios. After 16 hours, the apoptosis of NK cells was assessed. Assuming, that both cell populations are not able to grow inside the wells due to the lack of nutrients and space, we can use the following initial value problem to model the described phenomenon:

\begin{subequations} \label{NKinducedApoptosisByTregsSystemDataFit}
\begin{align}
    \dv{N}{t} &= -\theta_{N_E} N(t) - f\left(N,R\right) N(t)\,, \quad N(0) = N_E \,,\\
    \dv{R}{t} &= -\theta_{R_E} R(t)\,, \quad R(0) = ratio \cdot N_E\,,
\end{align}
\end{subequations}
where $N$ is the NK cell population, $R$ is the Treg population, $\theta_{N_E}$ is the rate of natural NK cell death \textit{in vitro}, $\theta_{R_E}$ is the rate of natural Treg death \textit{in vitro}, $N_E$ is the initial number of NK cells, $ratio$ is the ratio of Tregs to NK cells and $f\left(N,R\right)$ is the trophic function describing the killing of NK cells by Tregs.

From Figure 1a in \cite{mahnke2007depletion}, we notice that, \textit{in vitro}, the percentage of Treg reduction after 24 hours is approximately 18\%. Therefore, we have that $\theta_{R_E} \approx 0.1985$ day$^{-1}$. From Figure 5A in \cite{olkhanud2009breast}, we notice that, \textit{in vitro}, the percentage of NK cell reduction after 16 hours is approximately 39\%. Therefore, we have that $\theta_{N_E} \approx 0.7414$ day$^{-1}$ (see Appendix \ref{app:naturaldeath} for more on how to calculate the turnover rate \textit{in vitro}).

Since the authors of \cite{olkhanud2009breast} do not specify the initial number of NK cells that were put in the wells, we assume it to be $5 \cdot 10^9$, a number of the same order of
magnitude as the number of NK cells at homeostasis state we found in Appendix  \ref{SectionEquilibriumStates}. Choosing the following three trophic functions:

\begin{equation} \label{trophicFunctionsRelation2}
    f\left(N,R\right) = \gamma_N N^{\delta_N} \text{\quad or \quad} =\gamma_N \frac{N^{\delta_N}}{s_R R^{\delta_N} + N^{\delta_N}} \text{\quad or \quad}  =\gamma_N \frac{N}{\delta_N + N}\,,
\end{equation}
and following the same procedure we used for the NK-induced lysis of breast cancer cells in Section \ref{SectionLargeModelParameterEstimationTumourcells}, we get that the lytic curves predicted by  problem \eqref{NKinducedApoptosisByTregsSystemDataFit} are given in Figure \ref{fig:NKlysisByTreg}, whereas parameter values are given in Table \ref{NKinducedApoptosisByTregsParameterTable}. As we can see in Figure \ref{fig:NKlysisByTreg}, the power form and rational form make the greatest fit. Since the power form is simpler, we choose it over the rational form. For further simplicity, we round parameters $\gamma_N$ and $\delta_N$ in order to finally get $\gamma_N=10^{-6}\; \text{cell}^{-\delta_N} \cdot \text{day}^{-1}$ and $\delta_N = 0.5$.

\begin{table}
\centering
\caption{Parameter values admitted from data fitting problem \eqref{NKinducedApoptosisByTregsSystemDataFit} to data from \cite{shenouda2017exVivo}.
\label{NKinducedApoptosisByTregsParameterTable}}
\begin{tabular}{cc} 
\toprule
\multicolumn{2}{c}{\textbf{Power Form}}               \\
\midrule
$\gamma_N$ & $2.92131 \cdot 10^{-6}$                  \\
$\delta_N$ & 0.499502                                 \\ 
\midrule
\multicolumn{2}{c}{\textbf{Rational Form}}            \\
\midrule
$\gamma_N$ & $5.15405 \cdot 10^{10}$                  \\
$\delta_N$ & 0.478213                                 \\
$s_R$      & $2.93734 \cdot 10^{11}$                  \\ 
\midrule
\multicolumn{2}{c}{\textbf{Michaelis-Menten Form}}    \\
\midrule
$\gamma_N$ & 0.604742                                 \\
$\delta_N$ & $1.02378 \cdot 10^{10}$                  \\
\bottomrule
\end{tabular}
\end{table}

The rate of CD4$^+$ T cell-induced NK activation, $\kappa = 1.63 \cdot 10^{-11} \text{ cell}^{-1} \cdot \text{day}^{-1}$, is derived by assuming equation \eqref{LargeModeldN} at the zero-tumor equilibrium. When $\dv{N}{t}=0$, we have that $\sigma_N - \theta_N N_0  - \gamma_N R_0^{\delta_N} N_0   + \kappa H_0 N_0 = 0$ and by solving for $\kappa$ we derive its value.

The rate of NK cell death due to tumor interaction, $p=4.66 \cdot 10^{-12}\; \text{cell}^{-1} \cdot \text{day}^{-1}$, is derived by assuming equation \eqref{LargeModeldN} at the high-tumor equilibrium. When $\dv{N}{t}=0$, we have that $\sigma_N - \theta_N N_1 - pT_1 N_1 - \gamma_N R_1^{\delta_N} N_1   + \kappa H_1 N_1 = 0$ and by solving for $p$ we derive its value.

\subsection{The CD8$^+$ T cells} \label{SectionLargeModelParameterEstimationCD8cells}
The constant source of CD8$^+$ T cells, $\sigma_C = 3 \cdot 10^7 \; \text{cells} \cdot \text{day}^{-1}$, is taken from \cite{hellerstein}. In that study, the authors found that the mean value of the absolute proliferation of CD8$^+$ T cells is $5.9\; \text{cells} \cdot \text{day}^{-1} \cdot \mu$L$^{-1}$. Assuming that the average human has 5 liters of blood and converting the absolute proliferation of CD8$^+$ T cells to $\text{cells} \cdot \text{day}^{-1}$ and rounding the result up, yields the value of $\sigma_C$.

The rate of programmable CD8$^+$ T cell death, $\theta_C = 0.009$ day$^{-1}$, is found by assuming exponential decay of CD8$^+$ T cells and taking their half-life to be 77 days, as found in \cite{hellerstein}. Hence, $\theta_C = \tfrac{\ln{2}}{77 \text{ day}} \approx 0.009 \text{ day}^{-1}$.

The rate of CD8$^+$ T cell death due to tumor interaction, $q = 3.422 \cdot 10^{-10}$ cell$^{-1}$ $ \cdot $ day$^{-1}$, is borrowed from \cite{kuznetsov1994} in which the authors derived the value from mouse data and a general effector cell and cancer cell population.

The rate of Treg-induced CD8$^+$ T cell death, $\gamma_C = 10^{-6}\; \text{cell}^{-1} \cdot \text{day}^{-1}$, is an \textit{ad hoc} value that has been chosen to give reasonable biological results due to the lack of data regarding the death of CD8$^+$ T cells due to Tregs.

The rate of CD8$^+$ T cell activation due to NK lysed tumor cell debris, $r = 1.05 \cdot 10^{-10}\; \text{cells} \cdot \text{day}^{-1}$, is derived by assuming equation \eqref{LargeModeldC} at the high-tumor equilibrium. When $\dv{C}{t}=0$, we have that $\sigma_C  - \theta_C C_1 - qT_1C_1 - \gamma_C R_1C_1 + r N_1T_1 + \frac{j_C T_1}{k_C+T_1}C_1 + \frac{\eta_1 H_1}{\eta_2+H_1}C_1 = 0$ and by solving for $r$ we derive its value.

The rate of CD8$^+$ T cell recruitment due to cancer, $j_C = 1.245 \cdot 10^{-1}$ day$^{-1}$, and the breast cancer cell number for half-maximal CD8$^+$ T cell recruitment due to cancer, $k_C = 2.019 \cdot 10^7$ cells, are borrowed from \cite{kuznetsov1994} in which the authors derived the value from mouse data and a general effector cell and cancer cell population.

The rate of CD8$^+$ T cell recruitment due to CD4$^+$ T cells, $\eta_1 = 2.48$ day$^{-1}$, is derived by assuming equation \eqref{LargeModeldC} at the zero-tumor equilibrium. When $\dv{C}{t}=0$, we have that $\sigma_C  - \theta_C C_0 - \gamma_C R_0C_0 + \frac{\eta_1 H_0}{\eta_2+H_0}C_0 = 0$ and by solving for $\eta_1$ we derive its value.

The CD4$^+$ T cell number for half-maximal CD8$^+$ T cell recruitment due to CD4$^+$ T cells, $\eta_2 = 2.5036 \cdot 10^3$ cells, is an \textit{ad hoc} value that has been chosen to give reasonable biological results due to the lack of data regarding the activation of CD8$^+$ T cells by CD4$^+$ T cells.

\subsection{The CD4$^+$ T cells} \label{SectionLargeModelParameterEstimationCD4cells}
The constant source of CD4$^+$ T cells, $\sigma_H = 2.2 \cdot 10^7 $ cells · day$^{-1}$, is derived by assuming equation \eqref{LargeModeldH} at the zero-tumor equilibrium. When $\dv{H}{t}=0$, we have that $\sigma_H -\theta_H H_0 = 0$ and by solving for $\sigma_H$ we derive its value.

The rate of programmable CD4$^+$ T cell death, $\theta_H = 0.00797$ day$^{-1}$, is found by assuming exponential decay of CD4$^+$ T cells and taking their half-life to be 87 days, as found in \cite{hellerstein}. Hence, $\theta_H = \tfrac{\ln{2}}{87 \text{ day}} \approx 0.00797 \text{ day}^{-1}$.

The breast cancer cell number for half-maximal CD4$^+$ T cell recruitment due to breast cancer, $k_H = 2.5036 \cdot 10^3$ cells, is an \textit{ad hoc} value that has been chosen to give reasonable biological results due to the lack of data regarding the recruitment of CD4$^+$ T cells due to breast cancer cells.

The rate of CD4$^+$ T cell recruitment due to breast cancer, $j_H = 4.45 \cdot 10^{-12}$ cell$^{-1} \cdot$ day$^{-1}$, is derived by assuming equation \eqref{LargeModeldH} at the high-tumor equilibrium. When $\dv{H}{t}=0$, we have that $\sigma_H -\theta_H H_1 + \frac{j_H T_1}{k_H + T_1} B_1 H_1 - c_1 H_1 B_{T_1}  = 0$ and by solving for $j_H$ we derive its value.

The rate of differentiation of CD4$^+$ T cells to Tregs, $c_1 = 1.21 \cdot 10^{-10}$ cell$^{-1}$ · day$^{-1}$, is derived in Section \ref{SectionLargeModelParameterEstimationTregs}.

\subsection{The Tregs} \label{SectionLargeModelParameterEstimationTregs}
The constant source of Tregs, $\sigma_R = 9.24 \cdot 10^6 $ cells · day$^{-1}$, is derived by assuming equation \eqref{LargeModeldR} at the zero-tumor equilibrium. When $\dv{R}{t}=0$, we have that $\sigma_R -\theta_R R_0 = 0$ and by solving for $\sigma_R$ we derive its value.

The rate of programmable Treg death, $\theta_R = 0.03851$ day$^{-1}$, is found by assuming their half-life to be 18 days, as found in \cite{mabarrack2008}. Thus, assuming Tregs follow exponential decay we have that $ \theta_R = \frac{\ln2}{18 \text{ day}} \approx 0.03851 \text{ day}^{-1}. $

The rate of differentiation of CD4$^+$ T cells to Tregs, $c_1 = 1.21 \cdot 10^{-10}$ cell$^{-1}$ · day$^{-1}$, is derived by assuming equation \eqref{LargeModeldR} at the high-tumor equilibrium. When $\dv{R}{t}=0$, we have that $\sigma_R -\theta_R R_1 + c_1 H_1 B_{T_1} = 0$ and by solving for $c_1$ we derive its value.

\subsection{The B cells} \label{SectionLargeModelParameterEstimationBcells}
The constant source of non-tBreg B cells, $\sigma_B = 3.16 \cdot 10^7 $ cells · day$^{-1}$, is derived by assuming equation \eqref{LargeModeldB} at the zero-tumor equilibrium. When $\dv{B}{t}=0$, we have that $\sigma_B  -\theta_B B_0  = 0$ and by solving for $\sigma_B$ we derive its value.

The rate of programmable non-tBreg B cell death, $\theta_B = 0.0395$ day$^{-1}$, is derived from data taken from \cite{macallan2005}. In that study, the authors measured the half-life of the whole B cell population among 12 healthy donors aged between 19 and 85 years of age. Looking at the data from their Table 1, we have that the average B cell half-life in those 12 donors, with an average age of about 51.1 years, is approximately 17.56 days. Assuming exponential decay of B cells, we have that the rate of programmable B cell death is $ \theta_B = \frac{\ln2}{17.56 \text{ day}} \approx 0.0395 \text{ day}^{-1}. $

The rate of differentiation of B cells to tBregs, $c_2 = 1.7 \cdot 10^{-13}$ cell$^{-1}$·day$^{-1}$, is derived by assuming equation \eqref{LargeModeldB} at the high-tumor equilibrium. When $\dv{B}{t}=0$, we have that $\sigma_B  -\theta_B B_1 - c_2 T_1 B_1  = 0$ and by solving for $c_2$ we derive its value.

The rituximab-induced non-tBreg B cell inhibition coefficient, $\gamma_B = 20 \frac{\text{mL}^2}{\mu \text{g}^2 \; \cdot \; \text{day}}$, was found by running numerical simulations and choosing the value of $\gamma_B$, so the behavior of non-tBreg B cells to be similar to data from  \cite{tobinai1998feasibilityRituximabHalfLife} and \cite{cooper2004efficacyRituximab}. For a more in-depth discussion, see Section \ref{RituximabSimulationsSection}.

\subsection{The tBregs}  \label{SectionLargeModelParameterEstimationTBregs}
The rate of programmable tBreg death, $\theta_{B_T} = 0.039$ day$^{-1}$, is derived by assuming equation \eqref{LargeModeldBT} at the high-tumor equilibrium. When $\dv{B_T}{t}=0$, we have that $\theta_{B_T} B_{T_1} + c_2 T_1 B_1   = 0$ and by solving for $\theta_{B_T}$ we derive its value.

The rate of differentiation of B cells to tBregs, $c_2 = 1.7 \cdot 10^{-13}$ cell$^{-1}$·day$^{-1}$, was derived in Section \ref{SectionLargeModelParameterEstimationBcells}.

\subsection{The rituximab} \label{RituximabParamEstimationSection}
The rate of excretion of rituximab, $\theta_X$ = 0.033 day$^{-1}$, is taken from \cite{tobinai1998feasibilityRituximabHalfLife} and \cite{regazzi2005pharmacokineticRituximabHalfLife}. In \cite{tobinai1998feasibilityRituximabHalfLife}, the authors measured the half-life of 12 rituximab-treated  patients with relapsed CD20$^+$ B-cell lymphoma. Four of the patients received four weekly doses of 250mg/m$^2$ and eight of the patients received four weekly doses of 375mg/m$^2$. In total, the average rituximab half-life of both groups was 445.4 hours. In \cite{regazzi2005pharmacokineticRituximabHalfLife}, the authors measured the half-life of 22 patients with  follicular lymphoma in complete or partial remission, 14 patients with various autoimmune disorders, four patients with AL Amyloidosis and eight patients with relapsed follicular or mantle cell lymphoma. All patients received the standard dose of 375mg/m$^2$. Patients in the first two groups received four weekly doses, patients in the third group received eight weekly doses, whereas patients in the fourth group received a total of six doses with various schedules. No statistically significant difference was observed between the groups, with a total average half-life of about 3 weeks. Seeing as both studies agree on the half-life of rituximab being about 3 weeks and assuming exponential decay, we have that $ \theta_X = \frac{\ln2}{21 \text{ day}} \approx 0.033 \text{ day}^{-1}. $

The rituximab dose function, $v(t)$, is a function of time, and is determined as follows. The standard dosage of rituximab is 375mg/m$^2$ once a week for four weeks, as its clinical safety and efficacy has been established \cite{grillo2000rituximabStandardDose}. However, when inside the organism, rituximab is measured in $\mu$g/mL as we can notice in several studies, for example in \cite{tobinai1998feasibilityRituximabHalfLife} and \cite{regazzi2005pharmacokineticRituximabHalfLife}. Therefore, we convert the amount of rituximab received per dose from mg/m$^2$ to $\mu$g/mL.
We choose the body surface area to be equal to 1.7m$^2$ based on breast cancer patients' data provided by Table 3 in \cite{sacco2010average} and references there in. In  \cite{sacco2010average}, the authors focused on cancer patients who were already receiving some sort of treatment, but we do not notice a difference between them and patients who received no treatment based on values from other studies discussed in that particular article. Furthermore, assuming once again that the average human has 5 liters of blood, we have that

\begin{equation}
    375 \frac{\text{mg}}{\text{m}^2} \rightarrow 375 \cdot 1.7 \text{m}^2 \cdot \frac{1}{5 \text{L}}  \frac{\text{mg}}{\text{m}^2} = \frac{375 \cdot 1.7 }{5 \cdot 10^3 }  \frac{\text{mg}}{\text{mL}} = 127.5 \frac{\mu\text{g}}{\text{mL}}\,.
\end{equation}

Hence, we assume that every patient receives 127.5 $\mu$g/mL rituximab per dose.

In \cite{regazzi2005pharmacokineticRituximabHalfLife}, we see that the infusion time of rituximab is about 4 to 6 hours for the first infusion and 3 to 4 hours for subsequent infusions. Assuming that each infusion lasts 4 hours, we have that in order to model the total amount of rituximab entering the organism, the value of $v(t)$ needs to be equal to

\begin{equation}
    v(t) = \frac{127.5}{4} \frac{\mu\text{g}}{\text{mL} \cdot \text{hour}} = \frac{127.5}{4/24} \frac{\mu\text{g}}{\text{mL} \cdot \text{day}} = 765 \frac{\mu\text{g}}{\text{mL} \cdot \text{day}}\,,
\end{equation}
for 4 hours in order to simulate a full infusion. Therefore, in order to model a complete standard dose, starting at day 0, we have that

\begin{equation}
    v(t) = \begin{cases}
        765\,, \quad & t \in D \\
        0\,, \quad & \text{elsewhere}\,,
    \end{cases}
\end{equation}

where

\begin{equation}
        D = \{t\in\mathbb{R} : (0\le t \le 0.16) \cup (7\le t \le 7.16) \cup  (14\le t \le 14.16) \cup (21\le t \le 21.16) \}\,.
\end{equation}

To sum up we have Table \ref{tab:2}.

{\footnotesize
\begin{longtable}{p{0.7cm}p{5cm}p{2cm}p{2.2cm}p{3.5cm}}
\caption{A list of model parameters along with their description, value, units and source.}\\
\toprule
\textbf{Par.} & \textbf{Description} & \textbf{Value} & \textbf{Units} & \textbf{Source}\\
\hline
\endfirsthead
\multicolumn{5}{c}%
{\tablename\ \thetable\ -- \textit{Continued from previous page}} \\
\hline
\textbf{Par.} & \textbf{Description} & \textbf{Value} & \textbf{Units} & \textbf{Source} \\
\hline
\endhead
\hline \multicolumn{4}{r}{\textit{Continued on next page}} \\
\endfoot
\hline
\endlastfoot 
$a$                & Breast cancer growth rate                                                                & 0.17                   & day$^{-1}$                      & Data fitting from \cite{puchalapalli2016nsg}\\
$b$                & Inverse of breast cancer's carrying capacity                                             & $10^{-10}$             & cell$^{-1}$                     & Data fitting from \cite{puchalapalli2016nsg}\\
$\lambda_R$        & Treg-induced NK cell inhibition coefficient                                              & $10^{-8}$              & cell$^{-1}$                     & No data found\\
$c$                & Rate at which NK cells lyse breast cancer cells                                          & {[}11.2263, 19.6448]   & day$^{-1}$                      & Data fitting from \cite{shenouda2017exVivo}\\
$\delta$           & Hill coefficient measuring the steepness of the NK cell toxicity curve                   & {[}0.8249, 1.33332]    & -                               & Data fitting from \cite{shenouda2017exVivo}\\
$s_N$              & Value of $(\frac{N}{T})^\delta$ for half-maximal NK cell toxicity                          & {[}3.85119, 39.222]    & -                               & Data fitting from \cite{shenouda2017exVivo}\\
$d$                & Maximum rate at which CD8$^+$ T cells lyse cancer cells                                  & 1.7                    & day$^{-1}$                      & Borrowed from \cite{dePillis2013}\\
$l$                & Hill coefficient measuring the steepness of the CD8$^+$ T cell toxicity curve            & 1.7                    & -                               & Borrowed from \cite{dePillis2013}\\
$s_C$              & Value of $(\frac{C}{T})^l$ for half-maximal CD8$^+$ T cell toxicity                        & $3.5 \cdot 10^{-2}$    & -                               & Borrowed from \cite{dePillis2013}\\ 
\midrule                                                                                                      
$\sigma_N$         & Constant source of NK cells                                                              & $1.13 \cdot 10^8$      & cell $\cdot$ day$^{-1}$         & Estimated from \cite{zhang}\\
$\theta_N$         & Rate of programmable NK cell death                                                       & $0.06301$              & day$^{-1}$                      & Estimated from \cite{zhang}  \\
$p$                & Rate of NK cell death due to tumor interaction                                          & $4.66 \cdot 10^{-12}$   & cell$^{-1}$ · day$^{-1}$        & Estimated from homeostasis state \\
$\gamma_N$         & Rate of Treg-induced NK cell apoptosis                                                   & $10^{-6}$              & cell$^{-\delta_N} \cdot \text{ day}^{-1}$ & Data fitting from \cite{olkhanud2009breast} \\
$\delta_N $        & NK cell resistance to Treg-induced apoptosis coefficient                                 & $0.5$                  & -                               & Data fitting from \cite{olkhanud2009breast}\\
$\kappa$           & Rate of CD4$^+$-T-cell-induced NK activation                                             & $1.63 \cdot 10^{-11}$  & cell$^{-1}$ · day$^{-1}$        & Estimated from homeostasis state\\ 
\midrule                                                                                                      
$\sigma_C$         & Constant source of CD8$^+$ T cells                                                       & $3 \cdot 10^7$         & cell $\cdot$ day$^{-1}$         & Estimated from \cite{hellerstein}\\
$\theta_C$         & Rate of programmable CD8$^+$ T cell death                                                & 0.009                  & day$^{-1}$                      & Estimated from \cite{hellerstein}\\
$q$                & Rate of CD8$^+$ T cell death due to tumor interaction                                   & $3.422 \cdot 10^{-10}$  & cell$^{-1}$ · day$^{-1}$        & Borrowed from \cite{kuznetsov1994}\\
$\gamma_C$         & Rate of Treg-induced CD8$^+$ T cell death                                                & $10^{-6}$              & cell$^{-1}$ · day$^{-1}$        & No data found\\
$r$                & Rate of CD8$^+$ T cell activation due to NK lysed tumor cell debris                     & $1.05 \cdot 10^{-10}$   & cell$^{-1}$ · day$^{-1}$        & Estimated from homeostasis state\\
$j_C$              & Rate of CD8$^+$ T cell recruitment due to cancer                                         & $1.245 \cdot 10^{-1}$  & day$^{-1}$                      & Borrowed from \cite{kuznetsov1994}\\
$k_C$              & Breast cancer cell number for half-maximal CD8$^+$ T cell recruitment due to cancer      & $2.019 \cdot 10^7$     & cell                            & Borrowed from \cite{kuznetsov1994}\\
$\eta_1$           & Rate of CD8$^+$ T cell recruitment due to CD4$^+$ T cells                                & 2.48                   & day$^{-1}$                      & Estimated from homeostasis state\\
$\eta_2$           & CD4$^+$ T cell number for half-maximal CD8$^+$ T cell recruitment due to CD4$^+$ T cells & $2.5036 \cdot 10^3$    & cell                            & No data found\\ 
\midrule                                                                                                      
$\sigma_H$         & Constant source of CD4$^+$ T cells                                                       & $2.2 \cdot 10^7$       & cell · day$^{-1}$               & Estimated from homeostasis state\\
$\theta_H$         & Rate of programmable CD4$^+$ T cell death                                                & 0.00797                & day$^{-1}$                      & Estimated from \cite{hellerstein}\\
$j_H$              & Rate of CD4$^+$ T cell recruitment due to breast cancer                                  & $4.45 \cdot 10^{-12}$  & cell$^{-1} \cdot$ day$^{-1}$    & Estimated from homeostasis state\\
$k_H$              & Breast cancer cell number for half-maximal CD4$^+$ T cell recruitment                    & $2.5036 \cdot 10^3$    & cell                            & No data found\\
$c_1$              & Rate of differentiation of CD4$^+$ T cells to Tregs                                      & $1.21 \cdot 10^{-10}$  & cell$^{-1}$ · day$^{-1}$        & Estimated from homeostasis state\\ 
\midrule                                                                                                      
$\sigma_R$         & Constant source of Tregs                                                                 & $9.24 \cdot 10^6$      & cell · day$^{-1}$               & Estimated from homeostasis state\\
$\theta_R$         & Rate of differentiation of CD4$^+$ T cells to Tregs                                      & 0.03851                & day$^{-1}$                      & Estimated from \cite{mabarrack2008}\\ 
\midrule                                                                                                      
$\sigma_B$         & Constant source of non-tBreg B cells                                                     & $3.16 \cdot 10^7 $     & cell · day$^{-1}$               & Estimated from homeostasis state\\
$\theta_B$         & Rate of programmable non-tBreg B cell death                                              & 0.0395                 & day$^{-1}$                      & Estimated from \cite{macallan2005}\\
$c_2$              & Rate of differentiation of B cells to tBregs                                             & $1.7 \cdot 10^{-13}$   & cell$^{-1}$ · day$^{-1}$        & Estimated from homeostasis state\\ 
 $\gamma_B$ &  Rituximab-induced non-tBreg B cell inhibition coefficient & 20 & $\frac{\text{mL}^2}{\mu \text{g}^2 \; \cdot \; \text{day}}$        & Estimated from \cite{tobinai1998feasibilityRituximabHalfLife, cooper2004efficacyRituximab} \\
\midrule                                                                                                      
$\theta_{B_T}$     & Rate of programmable tBreg death                                                         & 0.039                  & day$^{-1}$                      & Estimated from homeostasis state\\ 
\midrule                                                                                     
$\theta_{X}$       & Rate of excretion of rituximab                                                           & 0.033                  & day$^{-1}$                      & Estimated from \cite{tobinai1998feasibilityRituximabHalfLife}, \cite{regazzi2005pharmacokineticRituximabHalfLife}
%\bottomrule
%\end{tabular}}
%\end{table}
\label{tab:2}  
\end{longtable}
}

\section{Numerical simulations and results}\label{numerics}
In this section, we numerically solve problem \eqref{LargeModelEquations}-\eqref{Model_ICs}  using Julia and the suite DifferentialEquations.jl \cite{rackauckas2017differentialequations}.
Before we begin, in order to get a better understanding of the breast tumor size, we convert the primary
tumor size classifications of the American Joint Committee on Cancer, found in Table 2 of \cite{giuliano2017breast}, from diameter (measured in mm) to total cell count. We present the results in Table \ref{BreastCancerStage}. We give a detailed explanation about this conversion in Appendix \ref{BreastCancerStageAppendix}.

\begin{table}
\centering
\caption{Classification of breast cancer size expressed in total breast cancer cell count. \label{BreastCancerStage}}
\begin{tabular}{cc} 
\toprule
\textbf{Category} & \textbf{Range (in Total Breast Cancer Cell Count)}  \\ 
\midrule
T1                & $[0, 2.30 \cdot 10^9]$                              \\
T1mi              & $[0, 2.88 \cdot 10^5]$                              \\
T1a               & $[2.88 \cdot 10^5, 3.59 \cdot 10^7]$                \\
T1b               & $[3.59 \cdot 10^7, 2.88 \cdot 10^8]$                \\
T1c               & $[2.88 \cdot 10^8, 2.30 \cdot 10^9]$                \\
T2                & $[2.30 \cdot 10^9, 3.59 \cdot 10^{10}]$             \\
T3                & $[3.59 \cdot 10^{10}, +\infty]$                     \\
\bottomrule
\end{tabular}
\end{table}

\subsection{Numerical simulations without rituximab} \label{NumericalSimulationsWithoutRituximab}
We begin by verifying whether our model yields a high-tumor equilibrium close to the high-tumor homeostasis values we calculated in  Appendix \ref{HighTumourEquilibrium}. As far as the zero-tumor equilibrium is concerned, there is no need to numerically verify its existence, since in our linear stability analysis in Appendix \ref{LinearStabilityLargeSection}, we calculated it analytically and, additionally, we derived some of the model's parameters using every coordinate of the zero-tumor homeostasis value we found in Appendix \ref{ZeroTumourEquilibriumSection}. In Figure \ref{fig:hightTumourHomeoastasis}, we see that with $E_1$ as the initial condition and parameter values as in Table \ref{tab:2}, with $c = 15$ day$^{-1}$, $s_N = 25$ and $\delta = 1$, the equilibrium of the model shows a slight decreased value of breast cancer cells when compared  to the biological homeostasis value determined in  Appendix \ref{HighTumourEquilibrium}. In particular, the number of breast cancer cells after 300 days is $9.97 \cdot 10^9$ cells, whereas all the other cells retain their initial value. This is to be expected, since we assumed breast cancer cells' high-tumor biological homeostasis value to be equal to that of an immunodeficient organism. Therefore, we verify that system \eqref{LargeModelEquations} exhibits biologically realistic results and we move on to studying the interactions of breast cancer and the immune system.

\begin{figure}[ht]
 \makebox[\textwidth][c]{\includegraphics[width=1\textwidth]{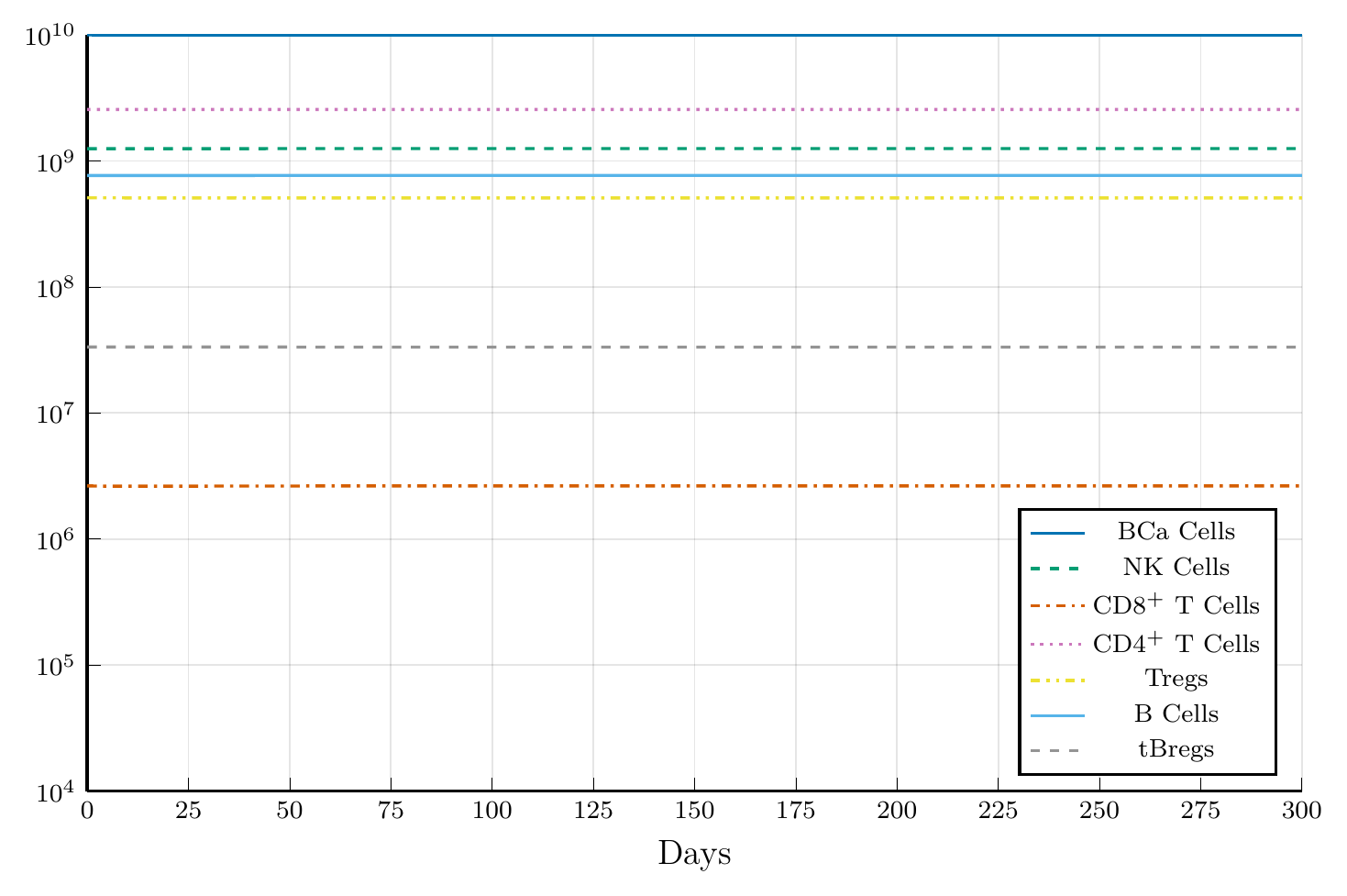}}%
  \caption{Initial condition equal to $E_1$ and parameter values as in Table \ref{tab:2}, with $c = 15$ day$^{-1}$, $s_N = 25$ and $\delta = 1$. }
  \label{fig:hightTumourHomeoastasis}
\end{figure}

\begin{figure}[ht]
 \makebox[\textwidth][c]{\includegraphics[width=1\textwidth]{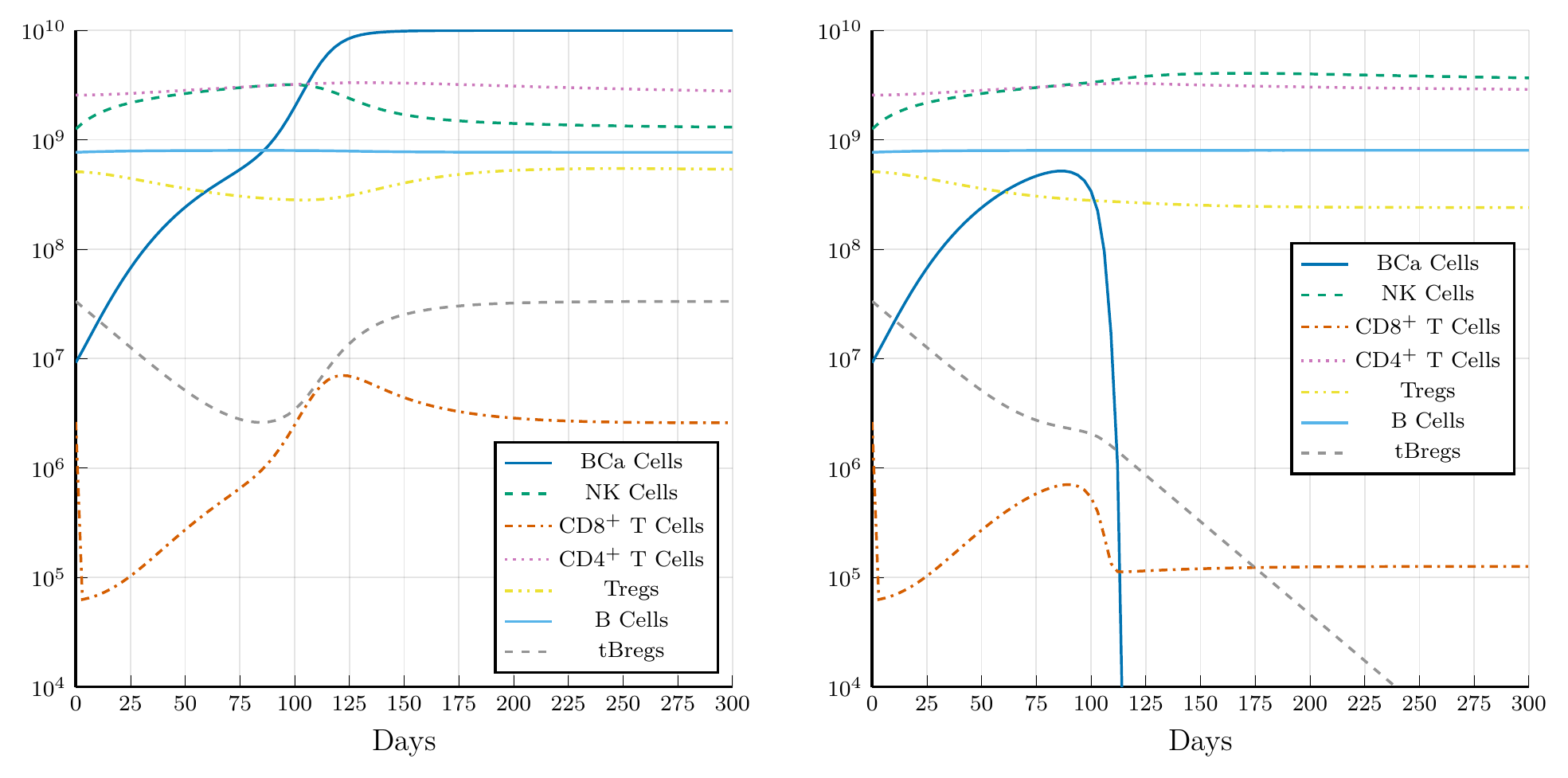}}%
  \caption{Initial condition of breast cancer cells is Left: $9.18 \cdot 10^6$ cells. Right: %Initial condition of breast cancer cells is
  $9.17 \cdot 10^6$ cells. All the other initial conditions are equal to $E_1$ and parameter values are as in Table \ref{tab:2}, with $c = 15$ day$^{-1}$, $s_N = 25$ and $\delta = 1$. }
  \label{fig:stabilityOfE1}
\end{figure}

Next, we numerically test the stability of the high-tumor equilibrium. Figure \ref{fig:stabilityOfE1} shows that an organism with the same parameter values as the simulation showcased in Figure \ref{fig:hightTumourHomeoastasis}, is not able to fight a relatively small, T1a-stage tumor, with a total cell population of $9.18 \cdot 10^6$ cells, while being able to kill any tumor lower than that. Evidently, the high-tumor equilibrium is stable.

\begin{figure}[t]
 \makebox[\textwidth][c]{\includegraphics[width=1\textwidth]{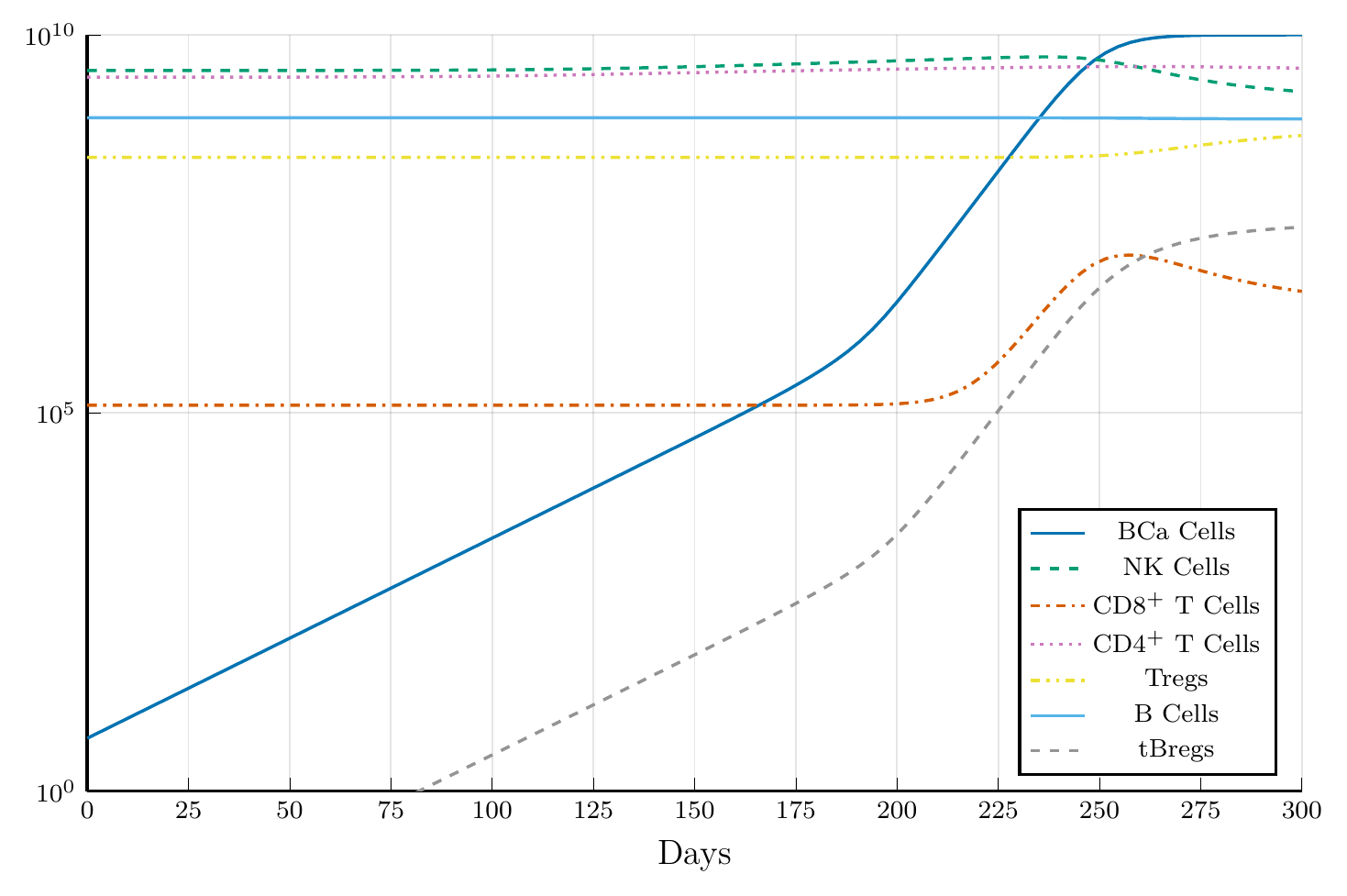}}%
  \caption{Initial condition of breast cancer cells is 5 cells. All the other initial conditions are equal to $E_0$ and parameter values, except for $d$ = 0.1, are as in Table \ref{tab:2}, with $c = 0.1$ day$^{-1}$, $s_N = 25$ and $\delta = 1$. }
  \label{fig:immuneBreakdown}
\end{figure}

We turn our attention to the zero-tumor equilibrium and its stability. As this equilibrium point's mathematical complexity is lower when compared to the complexity of the high-tumor equilibrium, we are able to analytically study its local stability. Linearization around the equilibrium shows that the zero-tumor equilibrium is locally stable for parameter values as in Table \ref{tab:2} (for more see Appendix \ref{LinearStabilityLargeSection}). In fact, out of the three parameters that represent NK cells' strength at lysing breast cancer cells, $c$, $s_N$ and $\delta$, only $c$ affects the local stability of the system and therefore causes the organism to either kill or succumb to the tumor when near the zero-tumor equilibrium. We note that $c$ is the only parameter of the aforementioned three, that does not directly relate to the ratio of NK to breast cancer cells, so it seems that in a healthy organism the rate at which NK cells lyse tumor cells is of greater importance than their ratio. That could also explain the case of cancer escaping immune surveillance and establishing itself, while only starting as a few cells. The same holds for the respective CD8$^+$ T cells parameters. For example, Figure \ref{fig:immuneBreakdown} shows a case of immune surveillance breakdown, with $c=0.1$ day$^{-1}$, $d=0.1$ day$^{-1}$ and an initial condition of 5 breast cancer cells while all the other cells are at their healthy homeostasis value. We see that after approximately 250 days, the tumor reaches its carrying capacity. An interesting observation is that in this scenario, a breast cancer tumor needs to be of around $10^3$ cells in order to generate tBregs and it does so at around 80 days after its formation.

\begin{figure}[t]
 \makebox[\textwidth][c]{\includegraphics[width=1\textwidth]{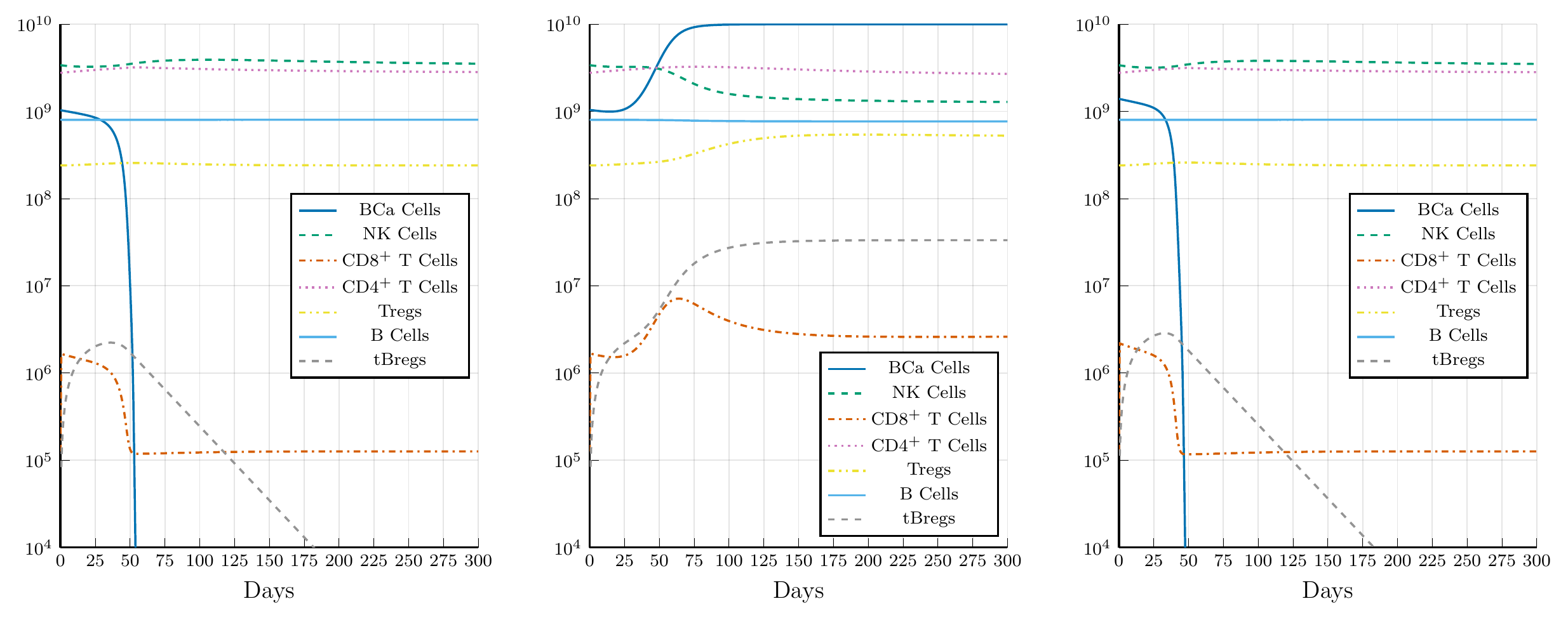}}%
  \caption{Left: Initial condition of breast cancer cells is $1.03 \cdot 10^9$ cells and $\delta =1$. Center: Initial condition of breast cancer cells is $1.04 \cdot 10^9$ cells and $\delta =1$. Right: Initial condition of breast cancer cells is $1.39 \cdot 10^9$ cells and $\delta =1.3$. All the other initial conditions are equal to $E_0$ and parameter values are as in Table \ref{tab:2}, with $c = 15$ day$^{-1}$ and $s_N = 25$. }
  \label{fig:stabilityOfE0}
\end{figure}

Being interested in whether $s_N$ and $\delta$ play a bigger role in tumor elimination as we move further away from the zero-tumor equilibrium, we run simulations to find what is the biggest tumor a healthy organism can beat. In Figure \ref{fig:stabilityOfE0}, we see that a healthy organism with parameter values as in Table \ref{tab:2} with $c=15$ day$^{-1}$, $s_N = 25$ and $\delta =1$, can kill T1c-stage tumors of around $1.03 \cdot 10^9$ cells, while unable to kill tumors larger than that. We also see that by increasing the value of $\delta$ to the maximum value found in our data fitting in Section \ref{SectionLargeModelParameterEstimationTumourcells}, the immune system is capable of beating tumors larger than the aforementioned size, while the same also holds for the case in which we decrease the parameter $s_N$ to its lowest, that is $s_N = 3.85119$. It is evident from the form of the functional response regarding the NK lysing of tumor cells, $-c \frac{N^\delta}{s_N T^\delta + N^\delta} $, that an increase in $c$ and a decrease in $s_N$, or in other words, an increase in the maximum rate at which NK cells lyse cancer cells and a decrease in the value of $(\frac{N}{T})^\delta$ for half-maximal NK toxicity, respectively, benefits the organism. Nevertheless, things are a bit more complicated as far as $\delta$ is concerned. In contrast to the results shown in Figure \ref{fig:stabilityOfE0}, Table \ref{dBad} shows the final number of cancer cells after 300 days in a simulation with initial conditions as in the high-tumor homeostasis values save for breast cancer cells and NK cells, all while allowing $\delta$ to take values outside of our data fitting results in order to illustrate how the tumor gets larger as $\delta$ increases. In order to make sense of these seemingly contradicting results, we take a closer look at the functional response term $-c \frac{N^\delta}{s_N T^\delta + N^\delta} $.

\begin{table}[b]
\centering
\caption{In all simulations the initial value of breast cancer and NK cells is $ 1.04 \cdot 10^9 $ cells and $ 5 \cdot 10^8 $ cells, respectively, with the other cells as in $E_1$. Parameter values are as in Table \ref{tab:2}, with $c=19$ day$^{-1}$ and $s_N =4$.}
\label{dBad}
\noindent\makebox[\textwidth]{%So it can be actually centred and not overfull. 
\begin{tabular}{lcccccc} 
\toprule
Value of $\delta$  & $ 0.0002$          & $ 0.002$           & $ 0.02$            & $ 0.2$             & $ 1$               & $ 2$                \\ 
BCa cells after 300 days & $8.182 \cdot 10^9 $ & $8.190 \cdot 10^9 $ & $8.259 \cdot 10^9 $ & $8.810 \cdot 10^9 $ & $9.777 \cdot 10^9 $ & $9.970 \cdot 10^9 $  \\
\bottomrule
\end{tabular}}
\end{table}

\begin{figure}[t]
 \makebox[\textwidth][c]{\includegraphics[width=1.0\textwidth]{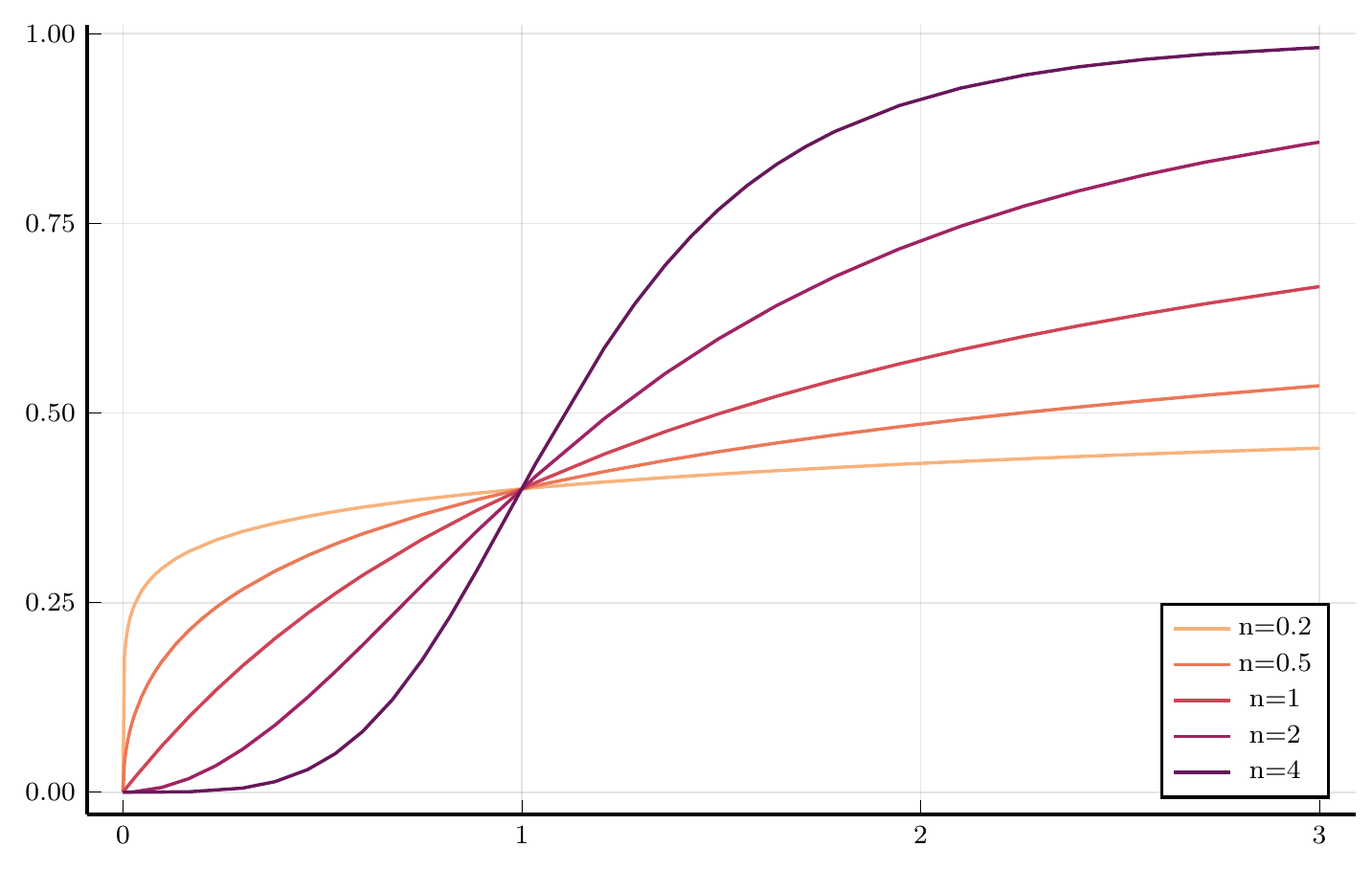}}%
  \caption{Plot of the Hill function for different values of the Hill coefficient, $n$. }
  \label{fig:HillPlot}
\end{figure}

\begin{figure}[ht]
 \makebox[\textwidth][c]{\includegraphics[width=1.0\textwidth]{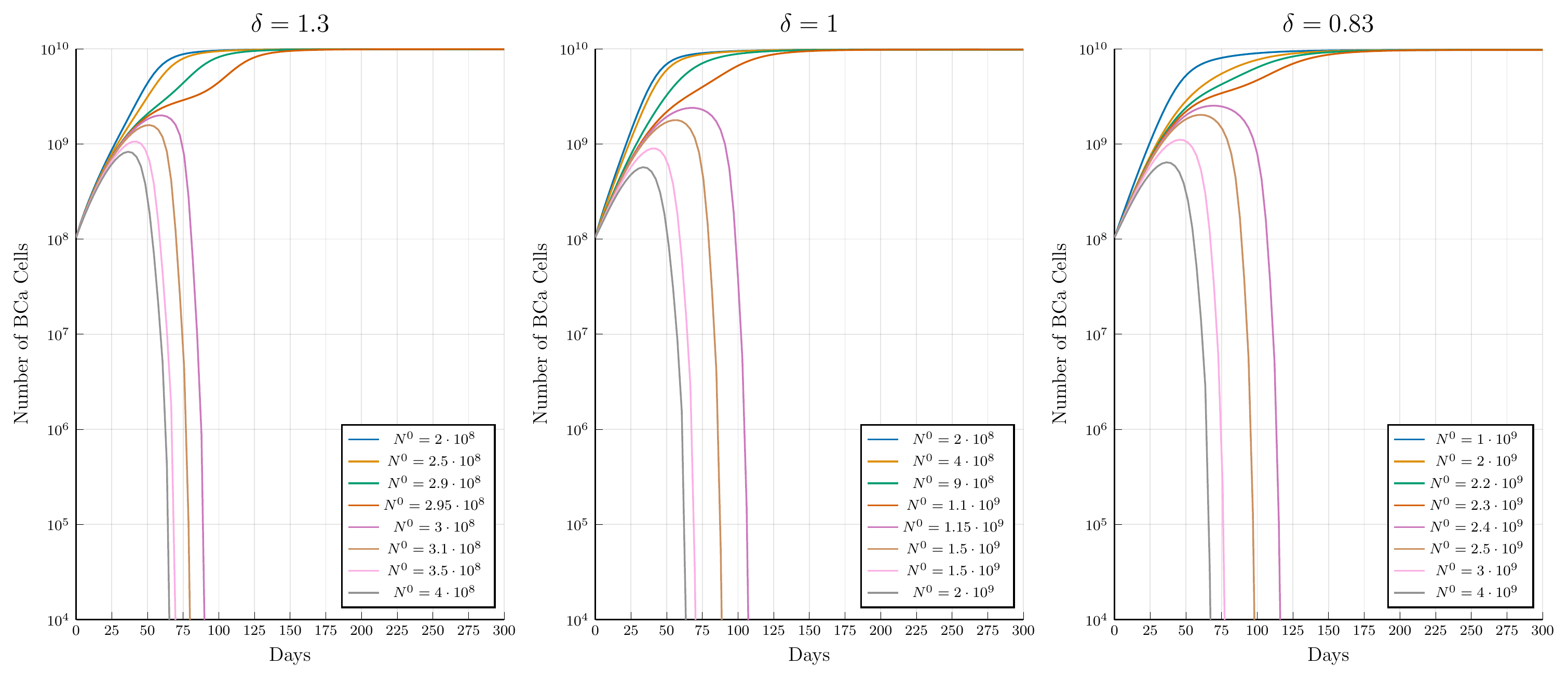}}%
  \caption{In all simulations the initial value of breast cancer cells is $1.04 \cdot 10^8$ cells and all the other cells are equal to their high-tumor homeostasis values. Parameters are as in Table \ref{tab:2} with $c=19$ day$^{-1}$ and $s_N = 4$.}
  \label{fig:changesInDelta}
\end{figure}

The functional response in question is a Hill function of the ratio of NK to breast cancer cells, with $\delta$ being its Hill coefficient. That means that values of $\delta$ smaller than 1 give to the curve of the breast cancer lysis by NK cells a hyperbolic form, whereas values of $\delta$ greater than 1 give to the curve a sigmoid form. This phenomenon is also present in our data fitting experiments in Figure \ref{fig:tumorLysisByNK}. Furthermore, as can be seen in Figure \ref{fig:HillPlot}, the smaller the Hill coefficient is, in this case $\delta$, the slower the lysis percent increases the more the NK to breast cancer cell ratio increases. On the other hand, the  larger  $\delta$ is, the faster the lysis percent increases when the ratio increases while near the curve's inflection point. In Figure \ref{fig:changesInDelta}, we showcase how increasing the initial value of NK cells affects the growing of breast cancer depending on the value of $\delta$.
The parameter $\delta$, as well as $c$ and $s_N$ in addition to breast-cancer-type-specific as we showed earlier, are also patient-specific, as found in \cite{dePillis2005validated} for the respective CD8$^+$ T cells case, and could theoretically be measured. Since one way of increasing the number of NK cells, and in turn the ratio, in real life could be by NK adaptive immunotherapy, measuring in advance the value of each patient's $\delta$, as well their total number of NK and breast cancer cells could be a valuable indication to whether an NK adoptive immunotherapy would have the intended results.

\begin{figure}[t]
 \makebox[\textwidth][c]{\includegraphics[width=1.0\textwidth]{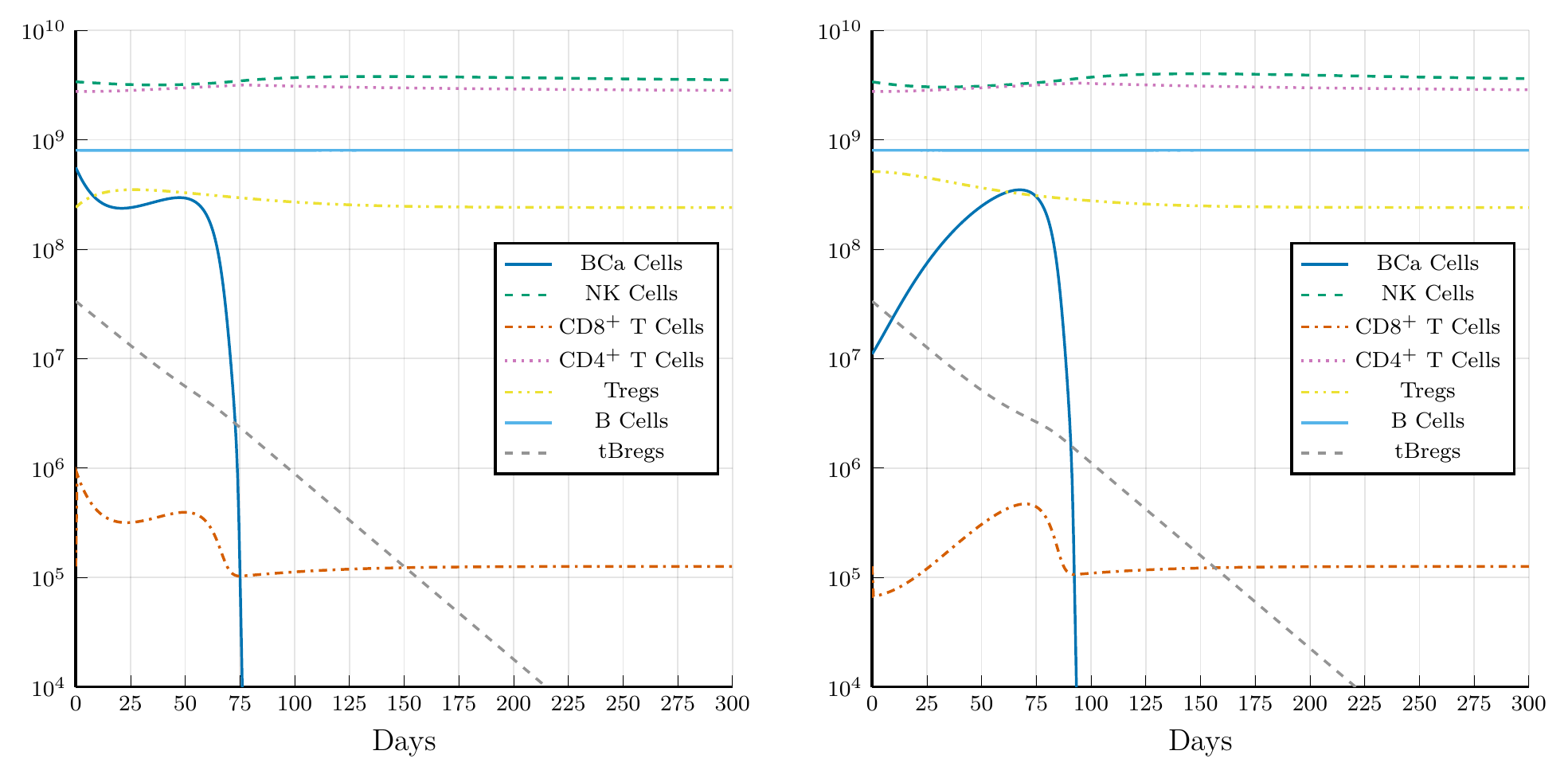}}%
  \caption{Left: Initial condition of breast cancer cells is $5.58 \cdot 10^8$ cells and of tBregs as their high-tumor homeostasis value. Right: Initial condition of breast cancer cells is $1.10 \cdot 10^7$ cells and of Tregs and tBregs as their high-tumor homeostasis value. In both simulations the initial value of all of the other cells is as their zero-tumor homeostasis values. Parameters are is in Table \ref{tab:2}, with $c=15$ day$^{-1}$, $s_N = 25$ and $\delta = 1$.}
  \label{fig:tBregInteraction}
\end{figure}

We continue our analysis by examining the interactions between breast cancer cells and tBregs. We saw earlier that a healthy organism could beat tumors as large as around $1.03 \cdot 10^9$ cells, which is a T1c-stage tumor. In the case, however, that the tumor has generated as many tBregs as their high-tumor equilibrium value, that number goes down to $5.58 \cdot 10^8$ cells as can be seen in Figure \ref{fig:tBregInteraction}, which is a T1b-stage tumor. When tBregs have led to the proliferation of Tregs that number goes even lower to $1.10 \cdot 10^7$ cells, which is a T1a-stage tumor. Clearly, another reason of breast cancer being able to establish itself is the result of the existence of regulatory cells. Our simulations show that the fewer they are the more likely it is for an organism to kill the tumor. The authors of \cite{olkhanud2011tumor} suggested the anti-CD20 antibody rituximab, which would deplete the B cell population and therefore tBregs - potentially stopping the proliferation of tBregs and in turn Tregs, as a possible therapy for breast cancer. We study the effects of rituximab in Section \ref{RituximabSimulationsSection}.

\subsection{Numerical simulations with rituximab} \label{RituximabSimulationsSection}

\begin{figure}
 \makebox[\textwidth][c]{\includegraphics[width=1.0\textwidth]{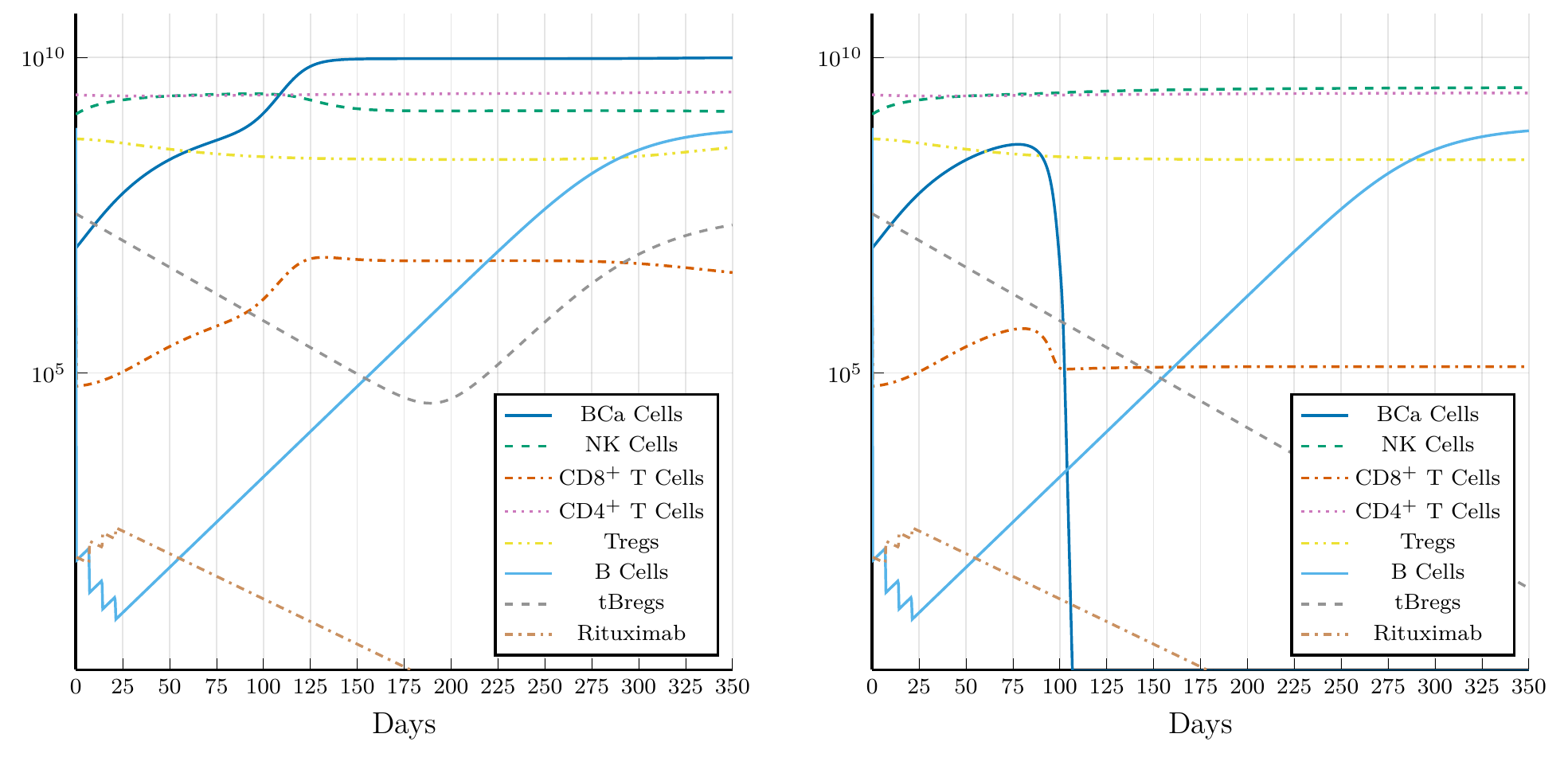}}%
  \caption{Initial condition of breast cancer cells is Left: $9.55 \cdot 10^6$ cells and Right: $9.54 \cdot 10^6$ cells. In both simulations the initial condition of all the other cells is at their high-tumor homeostasis value. Parameters are is in Table \ref{tab:2}, with $c=15$ day$^{-1}$, $s_N = 25$ and $\delta = 1$. Four weekly doses of 375 mg/m$^2$ of rituximab are administrated.}
  \label{fig:rituximabInitialBCa}
\end{figure}

We firstly check whether our model yields reasonable results in response to treatment with rituximab, based on data from clinical studies. In \cite{tobinai1998feasibilityRituximabHalfLife}, 12 patients with relapsed  CD20$^+$ B-cell lymphoma were treated with rituximab, as we already discussed in Section \ref{RituximabParamEstimationSection}. In all but one patients, B cells in the peripheral blood decreased to between 0 and 2\% of the total lymphocyte population, within two days after the first infusion. The remaining patient also exhibited reduced B cell levels, but only after the final infusion, that is four weeks after the first dose. In the four months that the patients' B cell levels were monitored, their B cell population did not recover. In \cite{cooper2004efficacyRituximab}, 57 patients with immune thrombocytopenic purpura were treated with the standard dose of rituximab, that is four weekly doses of 375 mg/m$^2$. As we can see in Figure 5 in \cite{cooper2004efficacyRituximab}, B cell levels started decreasing after the first infusion with rituximab and got depleted approximately five weeks later, after which B cells slowly increased until they regained their original population at around 51 weeks. As we can see in Figure \ref{fig:rituximabInitialBCa}, in our simulations, the non-tBreg B cell population rapidly decreases, just like in \cite{tobinai1998feasibilityRituximabHalfLife}, and reaches a population of around 10 cells, thus being almost depleted, at around the 25th day, that is around 10 days later when compared to the patients from \cite{cooper2004efficacyRituximab}. In our simulations, 120 days after the first dose, the non-tBreg B cell population is around 10$^4$ cells, which is negligible compared to the total lymphocyte population, just like the data from the two clinical studies. Furthermore, after 350 days the non-tBreg B cell population is slowly reaching its original number, that is about 10$^9$ cells, with the same thing happening at around the 357th day in \cite{cooper2004efficacyRituximab}. Since the behavior of non-tBreg B cells in our model is in agreement with the two clinical studies, we verify the validity of our model in predicting B cell behavior in response to rituximab treatment and we therefore proceed to analyze the results with regard to tumor growth.

In Figure \ref{fig:rituximabInitialBCa}, we see that with the standard treatment of four weekly doses of 375 mg/m$^2$ rituximab, the highest number of initial breast cancer cells an organism can beat raised slightly to $9.54 \cdot 10^6$ cells, when compared to $9.17 \cdot 10^6$ cells for the case without rituximab treatment, with both of them being T1a-stage tumors. We notice that the difference between the two breast cancer cell populations is very small. Even though tBregs exhibit a greater decrease when compared to the no rituximab case, it is still not enough for the organism to fight a significantly larger-sized tumor. Furthermore, tBregs start increasing shortly after the tumor has reached its carrying capacity. Additionally, we notice a slight decrease in the number of non-Treg CD4$^+$ T cells, which is to be expected since B cells activate CD4$^+$ T cells. We also notice a decrease in the Treg population, which returns to normal levels after tBregs reached their equilibrium. With these in mind, it seems that the tumor-induced differentiation of B cells to tBregs seems to play a bigger role than B cells activating CD4$^+$ T cells, as their depletion helps the organism, if only slightly. Hence, B cells seem to play a pro-tumor role in breast cancer growth.

\begin{figure}
 \makebox[\textwidth][c]{\includegraphics[width=1.0\textwidth]{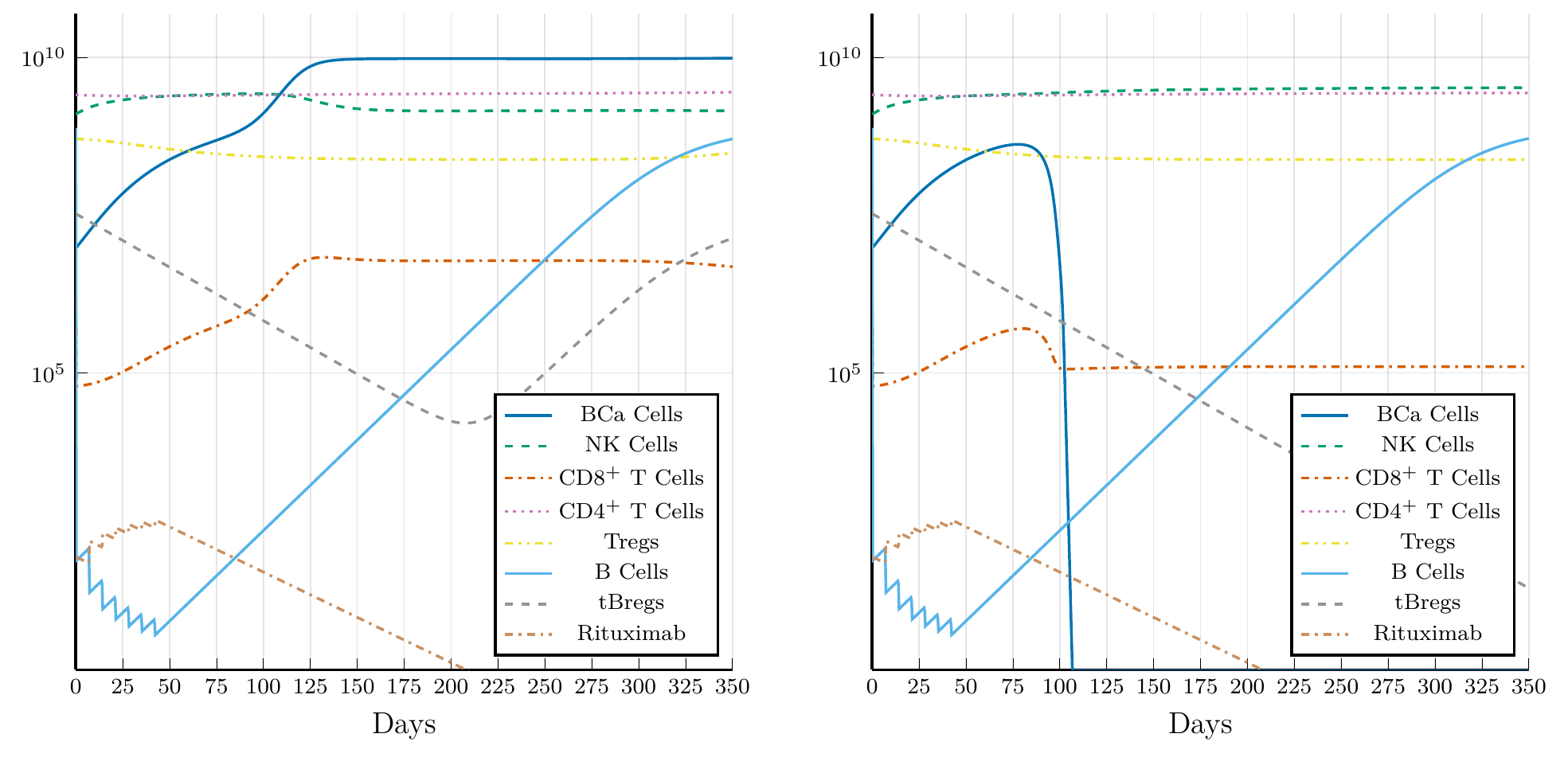}}%
  \caption{Initial condition of breast cancer cells is Left: $9.55 \cdot 10^6$ cells and Right: $9.54 \cdot 10^6$ cells. In both simulations the initial condition of all the other cells is at their high-tumor homeostasis value. Parameters are is in Table \ref{tab:2}, with $c=15$ day$^{-1}$, $s_N = 25$ and $\delta = 1$. Eight weekly doses of 375 mg/m$^2$ of rituximab are administrated.}
  \label{fig:rituximabInitialBCa8doses}
\end{figure}

We continue our analysis with trying out different experimental dosage schedules and quantities and evaluating their results. The administration of eight weekly doses of 375 mg/m$^2$ of rituximab, just like in \cite{regazzi2005pharmacokineticRituximabHalfLife}, does not change the maximum number of breast cancer cells an organism can beat, as we can see in Figure \ref{fig:rituximabInitialBCa8doses}. Even though that at four weekly doses tBregs decrease up to around the 188th day before they start increasing, at eight weekly doses the same things happens at around the 210th day. It is clear that more doses result in a further reduction of tBregs. However, this reduction happens too late, at a time when breast cancer has already reached its carrying capacity, thus making the organism unable to kill it.

\begin{figure}[t]
 \makebox[\textwidth][c]{\includegraphics[width=1.0\textwidth]{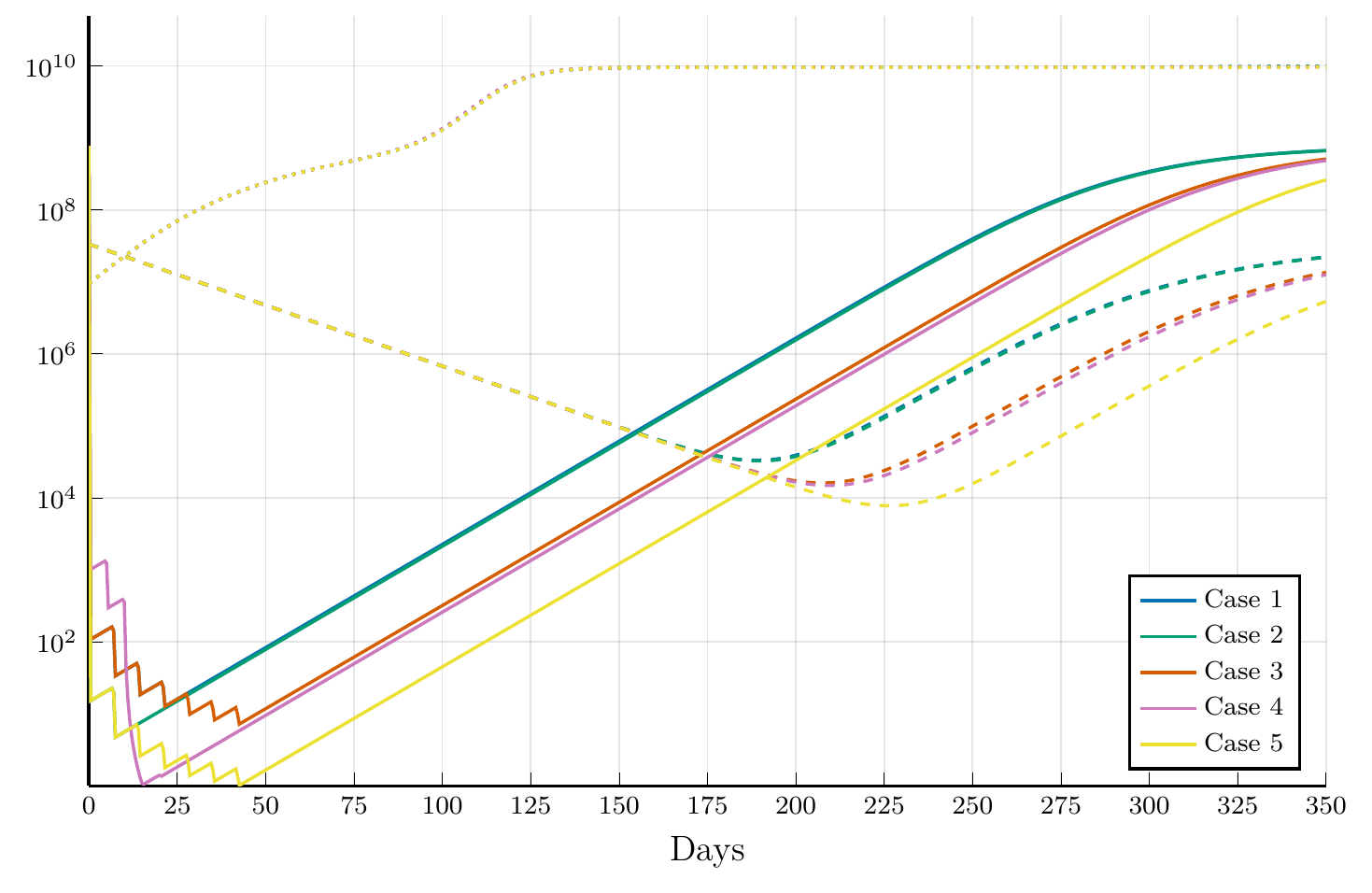}}%
  \caption{Initial condition of breast cancer cells is $9.55 \cdot 10^6$ cells, whereas all the other cells are at their high-tumor homeostasis values. Parameters are is in Table \ref{tab:2}, with $c=15$ day$^{-1}$, $s_N = 25$ and $\delta = 1$. A dotted line represents the tumor cell population, a dashed the tBregs and a solid the non-tBreg B cell population. Case 1: Four weekly doses of 375 mg/m$^2$. Case 2: Two weekly doses of 1 g/m$^2$. Case 3: Eight weekly doses of 375 mg/m$^2$. Case 4: Four doses of 122.549 mg/m$^2$ every five days. Case 5: Eight weekly doses of 1 g/m$^2$.}
  \label{fig:rituximabInitialBCa8doses3}
\end{figure}

In Figure \ref{fig:rituximabInitialBCa8doses3}, we compare five cases of different rituximab dosages. In the first case we consider a standard dose of four weekly doses of 375 mg/m$^2$, just like in Figure \ref{fig:rituximabInitialBCa}. For the second case, we increase the infused concentration of rituximab to 1 g/m$^2$ and decrease the number of infusions to two, just like in \cite{seitz2019highRituximab}, where the authors used the same dosage to treat patients with membranous glomerulonephritis. In the third case, we model eight weekly doses of 375 mg/m$^2$, like the case in Figure \ref{fig:rituximabInitialBCa8doses}. In the fourth case, we experiment with decreasing the quantity of rituximab infused to the patient to 122.549 mg/m$^2$, while also decreasing the dose schedule to be one infusion per five days. Finally, the fifth case is an extension of the third case, where we increase the dosage schedule to an 8-dose schedule and keep the quantity of rituximab to 1 g/m$^2$. In Figure \ref{fig:rituximabInitialBCa8doses3}, we notice that the treatment dosage that more successfully depletes non-tBreg B cells, while also being the superior at reducing tBregs is the fifth case. It is clear that no matter how successful each case is at depleting B cells, the effect it has on tBregs is to maintain the rate at which they decrease for a longer amount of time. In other words, the more successful in depleting the B cells a dosage is, the more time tBregs decrease for and they do so while maintaining the rate at which they decrease. However, that is not enough for the organism to beat the tumor, as that decrease of tBregs is not fast enough.

\subsection{Numerical sensitivity analysis} \label{NumericalSensitivityLargeSection}

\begin{figure}[t]
 \makebox[\textwidth][c]{\includegraphics[width=1\textwidth]{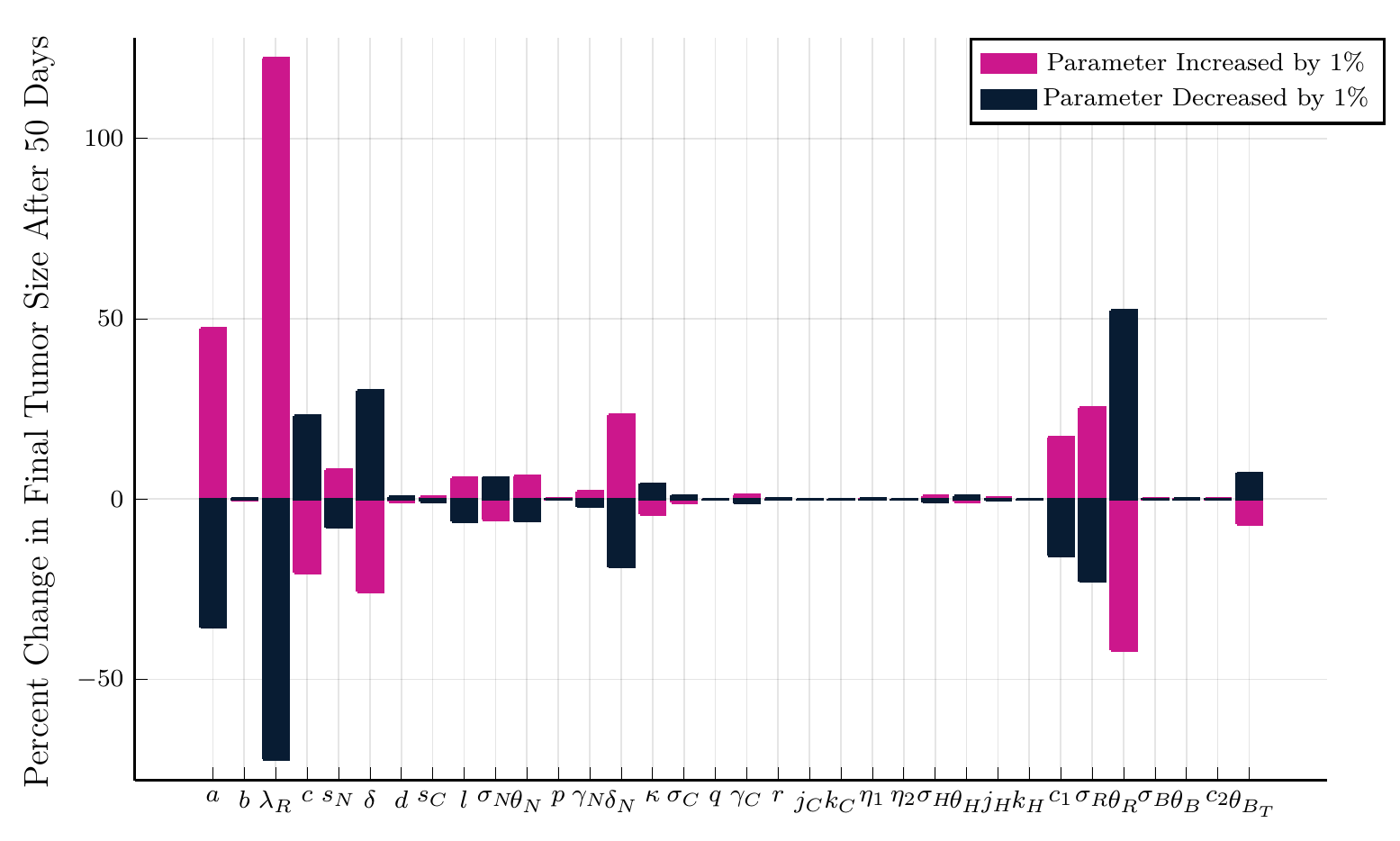}}%
  \caption{Depicted is the effect of a 1\% parameter change on final
tumor size after 50 days. Initial condition of breast cancer cells is $9.5 \cdot 10^6$ cells, whereas all the other cells are at their high-tumor homeostasis values. Parameters are is in Table \ref{tab:2}, with $c=15$ day$^{-1}$, $s_N = 25$ and $\delta = 1$.}
  \label{fig:LargeModelSensitiviy}
\end{figure}

In order to explore which parameters have the greatest effect on breast cancer-immune dynamics, we perform local sensitivity analysis on the model without rituximab interventions. The procedure has as follows. Firstly, we measure the final tumor size after 50 days with the initial condition of breast cancer cells being $9.5 \cdot 10^6$ cells, whereas all the other cells are at their high-tumor homeostasis values and parameter values as in Table \ref{tab:2}, with $c=15$ day$^{-1}$, $s_N = 25$ and $\delta = 1$. Next, we rerun the simulation, but this time we increase one parameter by 1\% and measure the percent change of the breast cancer cell population after 50 days when compared to our first simulation. Subsequently, we revert the parameter back to its original value and decrease it by 1\%, before rerunning the simulation and again measuring the percent change of the final tumor size when compared to our first simulation. After doing the same for all parameters, we get Figure \ref{fig:LargeModelSensitiviy}.

In Figure \ref{fig:LargeModelSensitiviy}, we notice that the parameters with the biggest impact on tumor growth, mainly concern five type of cells: breast cancer cells, NK cells, CD8$^+$ T cells, Tregs and tBregs. The parameter with the biggest impact on the system is the Treg-induced NK cell inhibition coefficient, $\lambda_R$, which is no surprise considering that it directly concerns three of the most important cells in our model: breast cancer cells, NK cells and Tregs. Finding a way to decrease the inhibition caused by Tregs to NK cells would greatly help the organism. It is also natural for the breast cancer growth rate, $a$, to play a big role in tumor growth. An interesting observation is that contrary to the parameters regarding NK cells, the parameters regarding CD8$^+$ T cells, show little to no sensitivity. Therefore, in the case studied in our sensitivity analysis, it is clear that NK cells are of greater importance when compared to CD8$^+$ T cells. Additionally, the three parameters we focused on in our numerical simulations, $c, s_N$ and $\delta$ also play a big role. Since, in our sensitivity analysis, the initial ratio of NK to breast cancer cells is by far greater  than one, increasing $\delta$, has an anti-tumor effect, while the opposite holds for decreasing $\delta$, further validating our claims in Section \ref{NumericalSimulationsWithoutRituximab}. Furthermore, the natural death rate of Tregs, $\theta_R$, their constant source, $\sigma_R$, as well as the rate of differentiation of CD4$^+$ T cells to Tregs, $c_1$, are also of importance. This implies that a drug such as sunitinib, which reduces the rate at which T cells differentiate into Tregs \cite{gu2010sunitinibDifferentiationTregs}, could prove useful in treating breast cancer. Finally, the reliance of our model on the natural death rate of tBregs, $\theta_{B_T}$, and the small reliance on parameters directly involving non-tBreg B cells, implies that a drug explicitly targeting tBregs, instead of implicitly targeting them through non-tBreg B cells, could be a better option.

\section{Conclusion and discussion}\label{conclusion}
In this study, we developed a model of nonlinear ordinary differential equations with the intent of exploring the various interactions between breast cancer and the immune system, with a focus on tBregs. Additionally, based on data fitting, we chose a Hill function with its variable being the ratio of NK cells to breast cancer cells, as the functional response which describes the way NK cells lyse breast cancer cells. 

Firstly, we validated the biological realism of our model by comparing its numerical solution with the two biologically realistic homeostasis values we derived.
Then, we found out that the largest tumor a healthy and a compromised organism could beat is a T1c-stage and a T1a-stage tumor, respectively.

Moreover, we gave the conditions under which an increase in the NK cell population, such as through the use of immunotherapy, could bear the intended results. 
These conditions revolved around the Hill coefficient of the functional response, which describes the way NK cells lyse breast cancer cells, 
as well as the ratio of NK cells to breast cancer cells. The Hill coefficient can be measured in a clinical setting for each particular patient though a chromium release assay, thus potentially making it a significant marker. 

Furthermore, we showed that when tBregs exist in an organism with the initial condition of all the other immune cells at the zero-tumor homeostasis state, the largest tumor the organism can beat goes down to a T1b-stage tumor and goes even lower to a T1c-stage tumor when tBregs have led to the proliferation of Tregs. Thus, we showed that tBregs need to be killed, in order for the tumor to be controlled.

We also performed simulations with the anti-CD20 antibody rituximab. After validating that the B cell decrease in our model mirrors that of clinical trials, we showed that with the standard rituximab dosage the size of the tumor an organism can beat increases, but does so only slightly. We additionally explored the behavior of our model with experimental rituximab dosages and found the same results as with the standard dose, as far as controlling the growth of breast cancer is concerned.

The aforementioned results, along with the reliance of breast cancer growth on tBregs and Tregs rather than directly on B cells, as was revealed by our sensitivity analysis, bore testament to the fact that attempts at controlling tBregs and Tregs could bring better results than targeting B cells.

\section*{Declaration of Competing Interest}
The authors declare that they have no known competing financial interests or personal relationships that could have appeared to influence the work reported in this paper.

\begin{appendices} \label{appendix}

\section{Homeostasis states} \label{SectionEquilibriumStates}
Here we derive two biological realistic homeostasis states for a zero-tumor condition and a high-tumor condition. Using the values of these two homeostasis states, we verify that our model yields biologically relevant results, as well as determine some of the model parameters.

\subsection{Zero-tumor homeostasis values} \label{ZeroTumourEquilibriumSection}
Naturally, $T_0 = 0$ in the zero-tumor homeostasis state.

Approximately 4 to 29\% of circulating lymphocytes are NK cells \cite{hema}. The average number of lymphocytes per microliter is 1000 to 4800 cells \cite{abbas2014cellular}, and since the average human has an average of 5 liters of blood \cite{starr2012biology}, we have that the total population of lymphocytes in a human is $5 \cdot 10^9$ to $24 \cdot 10^9$ cells. Therefore, the total population of NK cells in blood is $2 \cdot 10^8$ to $6.96 \cdot 10^9$ cells. Taking the median value yields $N_0 = 3.38 \cdot 10^9 $ cells. 

For the CD8$^+$ T cell zero-tumor homeostasis  value, we take the value derived from \cite{dePillis2009} which is $2.526 \cdot 10^4 \text{ cells} \cdot \text{L}^{-1}$, thus multiplying by 5 liters which is the average blood volume in a human, yields $C_0 = 1.263 \cdot 10^5 \text{ cells}$. This value represents the total number of CD8$^+$ T cells specific for a particular tumor associated antigen in the case of melanoma. While in this study  we are interested in breast cancer, the authors of \cite{dePillis2009} noted that other antigens present a similar degree of CD8$^+$ T cell activation.

The percentage of the total population of CD4$^+$ T cells among circulating lymphocytes ranges from 50 to 60\% \cite{abbas2014cellular}. Based on the total number of lymphocytes we calculated above, we have that the total number of circulating CD4$^+$ T cells ranges from $2.5 \cdot 10^9$ to $1.44 \cdot 10^{10}$ cells, so we choose the intermediate value of $3 \cdot 10^9$ cells. However, we are interested in the non-Treg CD4$^+$ T cell population, thus subtracting the median of the Treg population (see the following paragraph for its derivation) from the chosen intermediate value, we get the population of non-Treg CD4$^+$ T cells in the zero-tumor homeostasis state to be $H_0 = 2.76 \cdot 10^9$ cells.

Tregs make up 5 to 10\% of the circulating CD4$^+$ T cell population \cite{pang2013frequency}, or in other words 2.5 to 6\% of the whole circulating lymphocyte population, which means that their total population is in the range of $1.25 \cdot 10^8$ to $1.44 \cdot 10^9$ cells. Choosing the value corresponding to  8\% of the circulating CD4$^+$ T cell population \cite{liyanage} as our zero-tumor homeostasis state value for Tregs, we get $R_0 = 2.4 \cdot 10^8$ cells.

B cells are approximately 3 to 21\% of circulating lymphocytes in a healthy organism \cite{hema}, ergo their total population is $1.5 \cdot 10^8$ to $5.04 \cdot 10^9$ cells. Taking the median of that range and rounding it, we get that $B_0 = 8 \cdot 10^ 8$ cells.

We assume that there are no tBregs in the absence of tumor, at least not a clinically detectable number of them, therefore $B_{T_0} = 0$.

To summarize, the zero-tumor homeostasis state is
\begin{equation}
        E_0 = \big( 0,  3.38 \cdot 10^9, 1.263 \cdot 10^5, 2.76 \cdot 10^9,  2.4 \cdot 10^8, 8 \cdot 10^8, 0 \big) \cdot \text{cells}\,.
\end{equation}

\subsection{High-tumor homeostasis values} \label{HighTumourEquilibrium}
We assume that the breast cancer cell population at the high-tumor homeostasis state is equal to our model's carrying capacity parameter, $1/b$. Based on our data fitting in Section \ref{SectionLargeModelParameterEstimationTumourcells}, we choose the value of $b$ to be $10^{-10}$ cell$^{-1}$, thus the breast cancer cell population at the high-tumor homeostasis state is $T_1 = 1/b = 10^{10}$ cells.

For the the NK cell population, we use the same reasoning as in \cite{dePillis2009}, where authors noticed that in \cite{meropol1998evaluation} the average circulating NK cell population in cancer patients before receiving daily doses of IL-2 was 250 cells per microliter (data was taken from Figure 1 in \cite{meropol1998evaluation}), therefore $N_1 = 1.25 \cdot 10^9$ cells.

For the CD8$^+$ T cell high-tumor homeostasis state, we again use the value derived in \cite{dePillis2009}, 5.268$\cdot 10^5$ cells $ \cdot $ L$^{-1}$, thus $C_1 = 2.634 \cdot 10^6$ cells. Just like the homeostasis value of CD8$^+$ T cells in the zero-tumor homeostasis state, this value was calculated with data from CD8$^+$ T cells activated from a melanoma-specific antigen, however as the authors of \cite{dePillis2009} state, we assume a similar amount of  CD8$^+$ T cells to get activated in other types of cancers too.

In \cite{madu2013pattern} the authors measured the amount of circulating  CD4$^+$ T cells, in 80 cancer patients, 36 of which were suffering from breast cancer, before and 12 days after starting chemotherapy. Before starting chemotherapy the average number of circulating  CD4$^+$ T cells in breast cancer patients was 613 cells per microliter or about 3.065$\cdot 10^9$ cells in total. Since we are interested only in non-Treg CD4$^+$ T cells, we subtract $R_1$ (the Treg high-tumor homeostasis value found in the following paragraph) from the total CD4$^+$ T cell number, which yields $H_1 = 2.55621 \cdot 10^9$ cells.

In \cite{liyanage} the authors measured and compared the prevalence of Tregs in the whole CD4$^+$ T cell population among 35 breast cancer patients, 30 pancreatic cancer patients and 35 healthy donors. In the case of breast cancer patients they found that the percentage of Tregs among circulating CD4$^+$ T cells was higher when compared to healthy donors. Specifically, 16.6\% versus 8.6\%, respectively. Based on the 16.6\% prevalence of Tregs in breast cancer patients and the fact that we calculated that the average number of circulating CD4$^+$ T cells is 3.065$\cdot 10^9$ cells, we have that $R_1 = 5.0879 \cdot 10^8$ cells.

In \cite{tsuda2018b} the authors identified the expressions of cell markers from blood samples of 27 breast cancer patients and 12 healthy donors and found that the percentage of B cell in each cohort was about 8.905\% and 11.51\%, respectively (data was taken from Table 2 of \cite{tsuda2018b}). Those numbers are within the normal range of B cells \cite{hema,abbas2014cellular}, thus we'll assume the total B cell population remains constant when compared between a healthy person and a cancer patient, with the only changes happening within the B cell sub-populations, as we can also see in \cite{tsuda2018b}, with memory B cells being the most expanded sub-population. Consequently, since we are interested in the non-tBreg B cell population, we subtract the tBreg population assumed in the following paragraph from $B_0$ to find $B_1 = 7.67 \cdot 10^8 $ cells.

As tBregs are newly discovered, data regarding them are scarce. For that reason, we are unable to find the average population of tBregs in a breast cancer patient. Thus, we observe that in \cite{murakami2019increased} the authors discover an increase in the percentage of the immunosuppressive cytokine IL-10 producing B regulatory cells expressing the CD19$^+$CD24$^{\text{hi}}$CD27$^+$ mark in patients with gastric cancer compared to healthy donors. In particular about 8.35\% versus about 5.65\% of the whole CD19$^+$ expressing B cell population, respectively (data taken from Table 1 of \cite{murakami2019increased}). Since B cells express the CD19$^+$ mark, we consider this to be the whole B cell population. Furthermore, the authors of \cite{murakami2019increased} found out that the B cells expressing the CD19$^+$CD24$^{\text{hi}}$CD27$^+$ mark are able to suppress the proliferation of autologous CD4$^+$ T cells, while also inhibiting their IFN-gamma production. This makes us believe that there could be a connection between tBregs and CD19$^+$CD24$^{\text{hi}}$CD27$^+$ B cells, thus we make the assumption that half of those B regulatory cells are tBregs. Hence, multiplying the average number of B cells with 8.35\%, dividing by 2 and rounding it we get $B_{T_1} = 3.34 \cdot 10^7 $ cells.

To summarize, the high-tumor homeostasis state is
\begin{equation}
        E_1 = \big( 10^{10},  1.25 \cdot 10^9, 2.634 \cdot 10^6, 2.55621 \cdot 10^9,  5.0879 \cdot 10^8, 7.67 \cdot 10^8, 3.34 \cdot 10^7 \big) \cdot \text{cells}\,.
\end{equation}

\section{Calculation of the \textit{in-vitro} natural death rate} \label{app:naturaldeath}
Assuming that cells inside a well are not able to grow due to lack of nutrients and space, but only die due to natural death, we have that their population can be modeled by the following initial value problem:

\begin{equation} \label{cellPopulationInVitro}
    \dv{K}{t} = -\theta_{K_E} K(t)\,, \quad K(0) = K_E\,,
\end{equation}
where $K$ is the cell population, $\theta_{K_E}$ is the rate of natural cell death \textit{in vitro} and $K_E$ is the initial number of cells.

Assuming that at time $t_F$ we count the cell population and find that the population has been reduced by $p_E$\% compared to the initial cell population $K_E$, we have that $K(t_F) = K_E- \frac{p_E}{100} K_E $. By solving initial value problem \eqref{cellPopulationInVitro}, we get $K(t) = K_E e^{-\theta_{K_E}t}$. Setting $t=t_F$ and solving for $\theta_{K_E}$ yields

\begin{equation} \label{naturalcellDeathInVitroRelation}
    \theta_{K_E} = \frac{1}{t_F} \ln{\frac{1}{1-p_E}}\,.
\end{equation}

Relation \eqref{naturalcellDeathInVitroRelation} allows us to find the \textit{in vitro} natural death rate of a cell population, by only knowing the time that has passed since the cells were first put inside the wells until their assessment, and the percentage of their reduction. 

\section{Linear stability analysis of the zero-tumor equilibrium} \label{LinearStabilityLargeSection}

Even though system \eqref{LargeModelEquations} is too complex to analytically find all of its equilibria, we can calculate the equilibrium in which the tumor is zero, in the absence of rituximab. Let
\begin{equation}
    E^* = \left( T^*, N^*, C^*, H^*, R^*, B^*, B_T^* \right) \,,
\end{equation}
be the zero-tumor equilibrium.
Assuming that all derivatives are equal to zero and additionally that $T^*=0$, then from equation \eqref{LargeModeldB} we have that $B_T^* = 0$. From equation \eqref{LargeModeldR}, we get $R^* = \frac{\sigma_R}{\theta_R}$. From equation \eqref{LargeModeldB}, we get $B^* = \frac{\sigma_B}{\theta_B}$. From equation \eqref{LargeModeldH}, we get $H^* = \frac{\sigma_H}{\theta_H}$. Replacing $R^*$ and $H^*$ to equations \eqref{LargeModeldN} and \eqref{LargeModeldC}, we get $N^* = \frac{\theta _H s_N}{-\kappa  \sigma _H+\theta _H \theta _N+\theta _H \gamma _N \sqrt{\frac{\sigma _R}{\theta _R}}}$ and $C^* =\frac{\sigma _C}{\theta _C+\frac{\gamma _C \sigma _R}{\theta _R}-\frac{\eta _1 \sigma _H}{\eta _2 \theta _H+\sigma _H}}$.
Therefore, the zero-tumor equilibrium is 

\begin{equation}
    E^* = \Bigg(0, \; \frac{\theta _H s_N}{\Lambda},  \; \frac{\sigma _C}{\theta _C+\frac{\gamma _C \sigma _R}{\theta _R}-\frac{\eta _1 \sigma _H}{\eta _2 \theta _H+\sigma _H}}, \frac{\sigma _H}{\theta _H}, \; \frac{\sigma _R}{\theta _R}, \; \frac{\sigma _B}{\theta _B}, \;0  \Bigg)\,.
\end{equation}
 
 The Jacobian matrix of system \eqref{LargeModelEquations} at the equilibrium point $E^*$ is
 
\begin{equation}
\mathbf{J}(E^*) = 
\begin{bmatrix}
 A_{11} & 0 & 0 & 0 & 0 & 0 & 0 \\
 -\frac{p \theta _H s_N}{\Lambda} & -\frac{\Lambda}{\theta_H} & 0 & \frac{\kappa  \theta _H s_N}{\Lambda} & -\frac{\theta _H \gamma _N s_N}{2 \sqrt{\frac{\sigma _R}{\theta _R}} \Lambda} & 0 & 0 \\
 A_{31} & 0 & A_{33} & A_{34} & \frac{\gamma _C \sigma _C}{A_{33}} & 0 & 0 \\
 \frac{\sigma _B j_H \sigma _H}{\theta _B \theta _H k_H} & 0 & 0 & -\theta _H & 0 & 0 & -\frac{c_1 \sigma _H}{\theta _H} \\
 0 & 0 & 0 & 0 & -\theta _R & 0 & \frac{c_1 \sigma _H}{\theta _H} \\
 -\frac{c_2 \sigma _B}{\theta _B} & 0 & 0 & 0 & 0 & -\theta _B & 0 \\
 \frac{c_2 \sigma _B}{\theta _B} & 0 & 0 & 0 & 0 & 0 & -\theta _{B_T} \\
\end{bmatrix}\,,
\end{equation}
where

\begin{align*}
    & \Lambda = -\kappa  \sigma _H+\theta _H \theta _N+\theta _H \gamma _N \sqrt{\frac{\sigma _R}{\theta _R}}\,, \\
    & A_{11} = a-c e^{-\frac{\lambda _R \sigma _R}{\theta _R}}-d\,, \\
   & A_{22} =  \frac{\kappa  \sigma _H}{\theta _H}-\theta _N+\gamma _N \left(-\sqrt{\frac{\sigma _R}{\theta _R}}\right)\,, \\
    & A_{31} = \frac{\sigma _C \theta _R \left(j_C-q k_C\right) \left(\eta _2 \theta _H+\sigma _H\right)}{\gamma _C k_C \sigma _R \left(\eta _2 \theta _H+\sigma _H\right)+k_C \theta _R \left(\sigma _H \left(\theta _C-\eta _1\right)+\eta _2 \theta _C \theta _H\right)}+\frac{r \theta _H s_N}{\Lambda}\,, \\
    & A_{33} = -\theta _C-\frac{\gamma _C \sigma _R}{\theta _R}+\frac{\eta _1 \sigma _H}{\eta _2 \theta _H+\sigma _H} \quad \text{and} \\
    & A_{34} = \frac{\eta _1 \eta _2 \sigma _C \theta _H^2 \theta _R}{\left(\eta _2 \theta _H+\sigma _H\right) \left(\gamma _C \sigma _R \left(\eta _2 \theta _H+\sigma _H\right)+\sigma _H \theta _R \left(\theta _C-\eta _1\right)+\eta _2 \theta _C \theta _H \theta _R\right)}\,.
\end{align*}

The eigenvalues of $\mathbf{J}(E^*)$ are
\begin{equation}
    \lambda_1 = A_{11}\,, \; \lambda_2 = -\theta_B\,, \; \lambda_3 = -\theta_{B_T}\,, \; \lambda_4 = -\theta_H\,, \; \lambda_5 = -\theta_R\,, \lambda_6 = A_{33} \quad \text{and} \quad \lambda_7 = -\frac{\Lambda}{\theta_H}\,.
\end{equation}

For parameter values as in Table \ref{tab:2}, we have that $\lambda_i < 0,$ for $i=1,2,\dots,7$, hence $E^*$ is locally asymptotically stable.

\section{Classification of breast cancer size expressed in total breast cancer cell count} \label{BreastCancerStageAppendix}
The conversion of tumor diameter (measured in mm) to total cancer cell count that we present in Table \ref{BreastCancerStage} is found, using the same method as in Section \ref{SectionLargeModelParameterEstimationTumourcells} (taken from \cite{dePillis2013}), where we assume that a spherical cancer cell has a diameter of approximately 15.15$\mu$m, as well as that cancer cells and tumors are spherical. Next, we find the volume range for each tumor classification by utilizing the diameter of the largest tumor dimension, found in Table 2 of \cite{giuliano2017breast}. Finally, we divide the tumor volume by the cancer cell volume we derived, to find the range of each classification expressed in total breast cancer cell number.

\end{appendices}

%%%%%%%%%%%%%%%%%%%%%%%%%%%%%%%%%%%%%%%%%%%%%%%
\bibliographystyle{abbrv}
\bibliography{Bib_BT_2021}\label{bibliography}
\addcontentsline{toc}{chapter}{Bibliography}
%%%%%%%%%%%%%%%%%%%%%%%%%%%%%%%%%%%%%%%%%%%%%%%

\end{document}